\documentclass[11pt, reqno]{article}
\usepackage{amsmath}
\usepackage{amsthm}
\usepackage{enumitem}
\usepackage{amssymb}
\usepackage{latexsym}
\usepackage{color}
\usepackage{enumitem}

\newtheorem{thm}{Theorem}[section]

\newtheorem{lem}{Lemma}[section]
\newtheorem{prop}{Proposition}[section]

\newtheorem{rem}{Remark}[section]

\newtheorem{assumption}{Assumption}[section]
\numberwithin{equation}{section}

\begin{document}
	\title{Partially observable semi-Markov zero-sum stopping game with first-passage risk probability criterion}
	
	\author{Arnab Bhabak\thanks{Dhirubhai Ambani University, Gandhinagar, Gujarat-382007, India, email: arnab\_bhabak@dau.ac.in} \; Bivakar Bose\thanks{Department of Mathematics, Indian Institute of Technology Guwahati,
			Guwahati-781039, India, email: b.bivakar@iitg.ac.in}  \;   Chandan Pal\thanks{Department of Mathematics, Indian Institute of Technology Guwahati,
			Guwahati-781039, India, email: cpal@iitg.ac.in}   \;  
		and  Subhamay Saha\thanks{
			Department of Mathematics, Indian Institute of Technology,
			Guwahati, Guwahati-781039 India, email: saha.subhamay@iitg.ac.in}}
	\date{}
	
	\maketitle \baselineskip20pt
	\parskip10pt
	\parindent.4in
	\begin{abstract}
		\noindent In this paper, we study a partially observable zero-sum semi-Markov game with both stopping and control under first-passage risk probability criterion. At each decision epoch, if the observable state falls outside a prescribed target set, the players may either stop the game or choose to continue playing by selecting their actions. However, if the observable state enters the target set, the game terminates immediately. We prove that, under mild conditions, the value of the game is given by the unique solution of a system of dynamic programming inequalities with bilateral constraints. The main steps involve converting the partially observable game to a fully observable one by introducing a sequence of filtering measures that capture the joint distribution of the unobservable state and the remaining profit level. Furthermore, we establish the existence of a saddle point equilibrium of the game.
	\end{abstract}
	
	
	\section{Introduction}
	Markov decision processes (MDPs) and stochastic games are widely used to model decision-making in uncertain and dynamic environments, see \cite{puterman2014markov}. Stochastic games with fully observable states have been studied extensively in the literature, see \cite{parthasarathy2020stochastic} and the references therein. However, in many real-world situations, it can be seen that the true state of the system may not be fully observable and the decision-makers have to rely on partial information obtained through observation. Motivated by a wide range of applications in artificial intelligence, machine replacement, maintenance, signal processing, medicine, and many other fields (see \cite{deep2023partially,feinberg2016partially}), partially observable stochastic control problems and games have become an active area of research. In addition to the partial observation, this paper considers the underlying state process to be semi-Markov which provides a significant generalization, as the sojourn times between consecutive decision epochs may follow an arbitrary probability distribution. Apart from the classical risk-sensitive expected discounted and average-cost criteria, the risk probability criterion has emerged as another important criterion for addressing risk in stochastic control problems. This criterion measures performance by the probability that the total discounted accumulated payoff does not exceed a prescribed threshold. Because of this practical relevance it has been widely applied in many areas, such as the financial engineering and management sciences, etc.; for instance see \cite{boda2004stochastic,ohtsubo2002optimal,rendleman1987assessing,wen2022optimal} and the references therein. In many stochastic control problems, the horizon over which costs or rewards are accumulated may be uncertain rather than fixed (as finite or infinite). For example, in equipment maintenance problems, one is only interested upto the failure of the equipment. Such situations naturally motivate the study of first-passage models, where the objective is to optimize the accumulated cost or reward until the first time the state process enters a prescribed target set. These models have found important applications in reliability theory, maintenance planning, quality control, and related areas; see \cite{huang2011,mamabolo2004maintenance}.   Further, our game model has the feature that in addition to selecting an action, players have the opportunity to stop the game at any decision epoch, provided that at that epoch the observable state lies outside a given target set of observable states. Once the observable state reaches the target set, the game ends automatically. 
	
	Partially observable discounted infinite-horizon zero-sum games have been studied in \cite{ghosh2004partialzero}. The risk-sensitive partially observable setting has been investigated for MDPs under both finite and infinite-horizon criteria in \cite{baauerle2017partially}, and has subsequently been extended to zero-sum games in \cite{bhabak2024partially}. In the partially observable framework, the risk probability criterion has been considered for discrete-time finite-horizon stochastic games in \cite{wei2025partially}.
	For semi-Markov models, the risk probability criterion has been studied in both fully observable and partially observable cases. Finite-horizon semi-Markov decision processes with the risk probability criterion was investigated in \cite{huang2013minimum}, and the corresponding infinite-horizon problem was considered in \cite{wen2025riskkybernetika}. The partially observable counterpart of the finite-horizon model in \cite{huang2013minimum} was studied in \cite{wen2024semimarkovpartially}. First-passage risk probability criterion with fully observable semi-Markov processes have been considered in \cite{bhabak2022semimarkovprobabilitycriteria,huang2015firstpassageriskprobabilitysemimarkov}. In the partially observable setting, first-passage risk probability optimality for discrete-time zero-sum games have been recently studied in \cite{wei2026firstpassageriskprobability}. Stopping games, in which players can both control the evolution of the system and decide when to stop the game, have attracted considerable attention in recent years across a variety of performance criteria and stochastic frameworks. Under the discounted cost criterion, such models have been investigated in the discrete-time, continuous-time, and semi-Markov frameworks in \cite{ghosh2007zero}, \cite{pal2020continuous}, and \cite{chen2025zero}, respectively. The corresponding risk-sensitive game problems have been studied in the discrete-time and semi-Markov frameworks in \cite{zhang2024discrete} and \cite{bose2026semi}, respectively. In the partially observable framework, discrete-time finite-horizon stopping games under the risk probability criterion have recently been studied in \cite{zhang2026partialstoppingdiscrete}. Thus, to the best of our knowledge, this is the first work on partially observable semi-Markov stopping games under the first-passage risk probability criterion. In a broad sense, our work extends the setting of \cite{bose2026semi} by incorporating partial observability and the first-passage risk probability criterion. It also generalizes the work in  \cite{zhang2026partialstoppingdiscrete} by considering semi-Markov state dynamics and a random horizon induced by the first-passage criterion. Furthermore, it extends the problem studied in \cite{wei2026firstpassageriskprobability} to the stopping games framework.

	In this paper, we consider a zero-sum stopping game with a first-passage risk probability criterion, in which the underlying state process is semi-Markov and only partially observable. In the proposed game model, both the observable and unobservable state spaces are Borel spaces and the reward and terminal reward functions are non-negative but not necessarily bounded from above. Since the cost criterion depends on the unobservable state component, we introduce a sequence of filtering measures that characterize the joint conditional distribution of the unobservable state and the remaining profit level given the observable history. Based on these measures, we define an extended state process, introduce a new semi-Markov kernel and reformulate the original partially observable game as an equivalent completely observable semi-Markov game. To obtain the value of the game, we formulate a system of dynamic programming inequalities with bilateral constraints and show that it is equivalent to a single functional equation. We then show that the solution of this functional equation, and hence the function satisfying the system of inequalities, coincides with the value of the game. To construct the solution, we cannot employ the backwards-iteration method used in the finite-horizon setting of \cite{zhang2026partialstoppingdiscrete}, as our model involves a random time horizon. Instead, we construct a suitable operator and obtain the solution as the limit of a monotone sequence generated by this operator. Using the measurable selection theorem and Fan’s min-max theorem, then we define optimal control for both players. Finally, by constructing appropriate stopping times, we show that the resulting control–stopping pair forms a saddle-point equilibrium. Unlike in \cite{ghosh2007zero,pal2020continuous,zhang2024discrete,zhang2026partialstoppingdiscrete}, the semi-Markov dynamics considered in our model prevent the direct use of Dynkin's formula. Consequently, the proofs rely heavily on martingale-type arguments. Moreover, the presence of the risk probability criterion, together with the filtering measures introduced to handle partial observability, leads to additional technical challenges. As a result, the analysis follows a different line of argument from that in \cite{bose2026semi}. Thus, the contribution of this paper lies not only in the proposed model but also in the analytical techniques developed to establish the main results.

	The remainder of the paper is organized as follows. Section $2$ introduces the game model and the filtering measures. In Section 3, we prove that the equivalent functional equation admits a unique solution, show that this solution characterize the value of the game, and establish the existence of a saddle-point equilibrium.

	\section{Game Model}
	The partially observable semi-Markov zero-sum first-passage stopping game model under study is described as follows:
	\[
	\{X, Y, A, B, \{A(x)\subset A,\; B(x)\subset B : x\in X\}, Q, Q_{0}, r, \zeta_1, \zeta_2\},
	\]
	where each entity carries the following meaning:
	\begin{itemize}
		\item $X$, $Y$ are state spaces assumed to be Borel spaces endowed with the Borel
		$\sigma$-algebras $\mathcal{B}(X)$ and $\mathcal{B}(Y)$. $X$ is the observable state component, whereas $Y$ cannot be observed by both players.
		\item $A$ and $B$ are the set of all actions for Player $1$ and Player $2$, respectively, each
		taken as a Borel space with corresponding Borel $\sigma$-algebras
		$\mathcal{B}(A)$ and $\mathcal{B}(B)$.
		\item For each $x\in X$, $A(x)\in \mathcal{B}(A)$ and
		$B(x)\in \mathcal{B}(B)$, specify the admissible action sets in state $x$
		for Player $1$ and Player $2$, respectively. We assume that both $A(x)$ and $B(x)$ are finite for every $x\in X.$
		\item Define
		\[
		K := \{(x,a,b): x\in X,\; a\in A(x),\; b\in B(x)\}.
		\] 
		Given any $(x,a,b)\in K$ and $y \in Y$, the semi-Markov kernel $Q$ is defined on
		$[0,\infty)\times X \times Y$. Thus, if $a\in A(x)$ and $b\in B(x)$ are the actions
		chosen in observable state $x$, for the underlying state pair $(x,y)$, then for any $\theta>0$ and
		$E\in \mathcal{B}(X\times Y)$,
		$Q(\theta, E \mid x,y,a,b)$ is the joint probability that the sojourn time
		in state $(x,y)$ will be less than or equal to $\theta$ and the next state pair belongs to the set $E$.
		Moreover, for each $(x,a,b)\in K$,  $y \in Y$ and $E\in \mathcal{B}(X\times Y)$,
		$Q(\cdot,E\mid x,y,a,b)$ is a non-decreasing right-continuous real-valued
		measurable function on $[0,\infty)$ such that
		\[
		Q(0,E\mid x,y,a,b)=0.
		\]
		\item The initial distribution $Q_{0}$ on $\mathcal{B}(Y)$ is assumed to be known to both the players.
		\item $D\in X$ is a given target set, which can be thought of as the set of all bad states of the system,
		with $D^c:=X\setminus D$ denoting the complement of $D$ with respect to $X$.
		\item  We denote by $r$ the non-negative running reward function defined on $K\times Y$.
		\item  For $k=1,2$, $\zeta_k$ is the non-negative terminal cost function on $X\times Y$.
	\end{itemize}
	
	We now describe the evolution of the game.
	Suppose that the system state is $(x_{0},y_{0})$ at the initial
	decision time $t_{0}:=0$ where $x_{0}$ is observable and $y_{0}$ cannot be observed
	by the players, and they have a common profit goal $\lambda_{0}\in \mathbb{R}:=(-\infty, +\infty)$ in mind--that is, until the end of the game, Player $2$ aims to keep the accumulated cost within $\lambda_{0}$. Then based on the observable state $x_{0}\in D^c$, Player $1$ and Player $2$ independently and simultaneously choose actions $a_{0}\in A(x_{0})$ and $b_{0}\in B(x_{0})$ respectively. As a consequence
	of this choice, the system remains in $(x_{0},y_{0})$ until
	time (say) $t_{1}$, at which point the system state changes to (say) $(x_{1},y_{1})$, according
	to the transition law $Q(.,.\vert x_{0},y_{0},a_{0},b_{0})$.
	At time $t_{1}$, a reward $R(x_0,y_0,a_0,b_0,\theta_1):=\int_{0}^{\theta_1}e^{-\alpha s}r(x_0,y_0,a_0,b_0)ds=\frac{r(x_0,y_0,a_0,b_0)}{\alpha}(1-e^{-\alpha \theta_1})$ is earned by Player $1$ from Player $2$, where $\theta_{1}:=t_{1}-t_{0}$ and thus the remaining
	profit goal at $x_{1}$ is updated accordingly to  $\lambda_1:=(\lambda_0-R(x_0,y_0,a_0,b_0,\theta_1))e^{\alpha \theta_1}$. If $x_{1}\in D$ the game stops. If instead $x_{1}\in D^c$, then based on the current observable state $x_{1}$, the jump time $t_{1}$, as well as the previous state $x_{0}$, previous actions $a_0$ and $b_0$, Player $1$ and Player $2$ again independently and simultaneously choose actions $a_{1}\in A(x_{1})$ and $b_{1}\in B(x_{1})$, respectively, and the same sequence of events occurs.

	Instead of taking further actions, if the players choose to stop at a decision epoch corresponding to the current observable state $x_{k}\in D^c$, for some $k\in \mathbb{N}$, the following scenarios will occur:\\[0.4em]
	(i) If Player~1 decides to stop the game --- or if both players decide to stop simultaneously
	--- then Player~1 earns the terminal reward $\zeta_1(x_k, y_k)$ from Player~2.\\[0.2em]
	(ii) If Player~2 alone decides to stop the game, then Player~1 earns the terminal reward
	$\zeta_2(x_k, y_k)$ from Player~2.
	
	The process evolves in this way, and so we get an observable history $h_{n}$ of the game up to the $n$-th decision
	epoch,  namely, $h_{n}=(x_{0},a_{0},b_{0},\theta_{1},\ldots ,x_{n-1},a_{n-1},b_{n-1},\theta_{n},x_{n}),$ where, $(x_{m},a_{m},b_{m})\in K$, $y_{m}\in Y$ for $m=0,1,\ldots n-1$, and $x_{n}\in X$. We
	denote by $\mathcal{H}_{n}$ the set of all observable histories $h_{n}$ of the system up to the $n$-th decision epoch
	and assume $\mathcal{H}_{n}$ is endowed with Borel $\sigma$-algebras.
	
	To specify decision rules at each decision epoch, we need to describe strategies. For that, we introduce some notations. Let $\mathbb{R}_+:=[0,\infty)$. Define $\Omega := (X \times Y \times A \times B \times \mathbb{R} \times \mathbb{R}_+)^\infty$ and let $\mathfrak{F}$ be the corresponding Borel $\sigma$-algebra.
	
	 For each $n\geq 0$ and any trajectory $w=(x_{0},y_{0},a_{0},b_{0},\lambda_{0}, \theta_{1},\ldots,x_{n},y_{n},a_n, b_n,\lambda_{n},\theta_{n+1},\ldots)\in \Omega$, define 
	\begin{eqnarray}\label{stochastic process 8}
		&& X_{n}(w)=x_{n},\hspace{.2cm} Y_{n}(w)=y_{n},\hspace{.2cm} A_{n}(w)=a_{n},\hspace{.2cm} B_{n}(w)=b_{n},\nonumber \\
		&&\Theta_{n}(w)=\theta_{n},\hspace{.2cm} T_{0}(w)=0,\hspace{.2cm} T_{n+1}(w)=t_{n+1}=\theta_{1}+\ldots+\theta_{n+1},\\&&\Lambda_0=\lambda_0,\quad \Lambda_n(w)=\lambda_n=(\lambda_{n-1}-R(x_{n-1},y_{n-1},a_{n-1},b_{n-1},\theta_n))e^{\alpha \theta_n},\nonumber
	\end{eqnarray}
	where $X_{n}$ denotes the observable state at $n$-th decision epoch, $Y_{n}$ denotes the unobservable state at $n$-th decision epoch. The $n$-th jump time and sojourn time at
	the $n$-th state are denoted by $T_{n}$ and $\Theta_{n+1}$, respectively. At the time $T_n$, the remaining budget and the actions chosen by the Player $1$ and Player $2$ are denoted by $\Lambda_n$, $A_n$ and $B_n$, respectively.  Furthermore, we define the natural filtration generated by the processes as:
	$$
	\mathcal{F}_0 :=\sigma\left(H_0\right),\quad	\mathcal{F}_n := \sigma\left(H_n\right),~\forall n \geq 1,
	$$ where $H_0:=X_0$ and $H_n:=(X_0,A_0,B_0,T_1,\ldots,X_{n-1}, A_{n-1}, B_{n-1}, T_{n},X_n),~ \forall n\geq 1$.

	\noindent\textbf{Stopping Time: } A random variable $\tau: \Omega\rightarrow \mathbb{N}\cup \{0\} \cup \{\infty\}$ is called a stopping time if
	$$\{\tau =n\}\in \mathcal{F}_{n}, \hspace{.2cm} \forall n\geq 0.$$
	We denote the set of all stopping times by $\mathcal{T}$.
	
	It is natural to consider the policies that only depend on the
	observed histories. So we define the following randomized history-dependent policy for the game, which specifies the action selection rules for the decision-makers.\\
	\textbf{Strategies:} For Player $1$ a randomized history dependent strategy $\pi^1:=\{\pi_{n}^1\}_{n\geq 0}$ is a
	sequence of stochastic kernels $\pi_{n}^{1}$ on $A$ given $\mathcal{H}_{n}$ satisfying $\pi_{n}^{1}(A(x_{n})\vert h_{n})=1$, $\forall h_{n}\in \mathcal{H}_{n}$. The family of all randomized history-dependent strategies for Player $1$ is denoted by $\Pi_{1}$. Similarly, we can define a randomized history-dependent strategy for Player $2$. Let
	$\Pi_{2}$ denotes the set of all randomized history-dependent strategies for Player $2$.
	In this particular game setup a strategy for Player $k$, is a pair $(\pi^k,\tau_k)$, where $\pi^k\in \Pi_{k}$ and $\tau_{k}\in \mathcal{T}$, $k=1,2$.
	
	For any  $x \in X$,  $y\in Y$, $\lambda \in \mathbb{R}$ and $\pi^k \in \Pi_k$,  $k=1,2$ by the well-known Tulcea's theorem \cite[Proposition 7.28]{Bertsekas1996stochastic}, there exists a unique probability measure $ \mathbb{P}_{x,y,\lambda}^{\pi^1,\pi^2}$ on the measurable space  $\left(\Omega, \mathfrak{F}\right)$ such that, for each $\theta \in \mathbb{R}_+,~ \mathbb{A} \in \mathcal{B}(A),~ \mathbb{B} \in \mathcal{B}(B),~ E\in \mathcal{B}(X\times Y),~ h_n\in \mathcal{H}_n$ and $n \geq 0$,
	\begin{equation}\label{property of probability measure 8}
		\begin{aligned}
			&\mathbb{P}_{x,y,\lambda}^{\pi^1,\pi^2}(X_0=x, Y_0=y, \Lambda_0=\lambda)=1\\& \mathbb{P}_{x,y,\lambda}^{\pi^1,\pi^2}\left(A_n \in \mathbb{A}, B_n \in \mathbb{B}\mid h_n\right)=\pi^1_n(\mathbb{A}\mid h_n)\pi^2_n(\mathbb{B}\mid h_n) \\
			& \mathbb{P}_{x,y,\lambda}^{\pi^1,\pi^2}\left(\Theta_{n+1}\leq \theta, (X_{n+1}, Y_{n+1})\in E \mid h_n,y_n, a_n,b_n\right)=Q\left(\theta, E \mid x_n,y_n, a_n, b_n\right)\\
			&\mathbb{P}_{x,y,\lambda}^{\pi^1,\pi^2}(\Lambda_{n+1}=(\lambda_{n}-R(x_{n},y_{n},a_{n},b_{n},\theta_{n+1}))e^{\alpha \theta_{n+1}}\mid h_n,y_n,a_n,b_n,\theta_{n+1})=1.
		\end{aligned}
	\end{equation}
	
	Let, $\mathbb{P}_{x,\lambda}^{\pi^1,\pi^2}(\cdot):=\int_{{Y\times \mathbb{R}}}\mathbb{P}_{x,y,\ell}^{\pi^1,\pi^2}(\cdot)Q_0(dy)\delta_{\lambda}(d\ell).$ Moreover,	let $\mathbb{E}_{x,\lambda}^{\pi^1,\pi^2}$ denote the expectation with respect to the probability measure $\mathbb{P}_{x,\lambda}^{\pi^1,\pi^2}$.
	
	In order to avoid the possibility of an infinite number of jumps within a finite time interval, we make the following standard assumption.
	\begin{assumption}\label{assumption 8.1}
		There exist constants $\delta>0$ and $\epsilon>0$ such that
		$$
		\sup_{(x,y,a,b) \in K} Q(\delta, X,Y \mid x,y, a,b) \leq 1-\epsilon
		$$
	\end{assumption}
	
	\begin{rem}\label{remark 2 finite jump}
		Assumption \ref{assumption 8.1} tells that the sojourn time at any state and under any action exceeds $\delta$ with probability at least $\epsilon$. If $T_{\infty}=\lim _{n \rightarrow \infty} T_n$, then it is well known that (see \cite[Proposition 2.1]{huang2011performance}), if Assumption \ref{assumption 8.1} holds then $\mathbb{P}_{x,y,\lambda}^{\pi^1,\pi^2}\left(T_{\infty}=\infty\right)=1$ for any $(x,y) \in X\times Y$, $\lambda \in \mathbb{R}$ and $\pi^k \in \Pi_k$, $k=1,2$.
	\end{rem}
	We now proceed to evaluate payoffs. A distinctive feature of this class of games is that the planning horizon is random, which arises because either player may stop the game at any decision epoch, or the system may enter the observable set of states $D$. The cost accumulates over the course of the game. Finally, as Player $1$ aims to exceed the prescribed profit threshold, while Player $2$ seeks to keep the accumulated cost within it, Player $1$ minimizes and Player $2$ maximizes the random-horizon first-passage risk-probability function defined below.
	
	\noindent {\bf Risk-probability criteria:} Fix a discounting factor $\alpha>0$ and initial profit level $\lambda \in \mathbb{R}$. Then,  for any pair of strategy $\left((\pi^1,\tau_1),(\pi^2,\tau_2)\right)$ $\in \left(\Pi_1\times \mathcal{T}\right)\times  \left(\Pi_2\times \mathcal{T}\right)$ and initial observable state $x \in X$, the risk-probability function is given by	
	\begin{equation}\label{cost criterion 8}
		\begin{aligned}
			&\mathbb{J}_\alpha\left(x,\lambda, \pi^1, \tau_1, \pi^2, \tau_2\right) \\
			&:=\int_Y\mathbb{P}_{x,y,\lambda}^{\pi^1, \pi^2}\bigg\{\sum_{n=0}^{(\tau_1 \wedge \tau_2\wedge \tau_D) - 1}
			e^{\!-\alpha T_n} R\left(X_n, Y_n, A_n, B_n, \Theta_{n+1}\right) \\
			&\hspace{1.2cm}+ \bigg(e^{\!-\alpha T_{\tau_1 \wedge \tau_2}}\bigg(\zeta_1\left(X_{\tau_1},Y_{\tau_1}\right) I_{\{\tau_1 \leq \tau_2\}}
			+ \zeta_2\left(X_{\tau_2},Y_{\tau_2}\right) I_{\{\tau_1 > \tau_2\}}\bigg)\bigg)I_{\{\tau_1\wedge \tau_2\leq \tau_D-1\}}\leq \lambda
			\bigg\}Q_0(dy),
		\end{aligned}
	\end{equation} where $\tau_D:=\inf\{n\geq 0: X_n\in D\}$ is the first passage time to the target set $D$.
	
	In \eqref{cost criterion 8} and what follows, we use the following convention.
	For any sequence $\left\{\xi_k\right\}_{k\geq 0}$ and $n,m\in \mathbb{N}$, 
	\begin{equation}\label{convention 8}
		\sum_{k=n}^m \xi_k:=0,~ \text{if}~m<n.
	\end{equation}
	
	\begin{rem}
		Using the convention in \eqref{convention 8}, from \eqref{cost criterion 8} it follows that, for all $(x,\lambda)\in D\times \mathbb{R}$ and $\left((\pi^1,\tau_1),(\pi^2,\tau_2)\right) \in \left(\Pi_1\times \mathcal{T}\right)\times  \left(\Pi_2\times \mathcal{T}\right)$, we have 
		$$
		\mathbb{J}_\alpha\left(x,\lambda, \pi^1, \tau_1, \pi^2, \tau_2\right)=I_{[0,\infty)}(\lambda). 
		$$ Hence, if the initial observable state $x\in D$, the first passage risk-probability function doesn't depend on the players' strategies and stopping times. Therefore, we shall henceforth restrict our analysis to initial observable states $x\in D^c$.
	\end{rem}
	For each $(x,\lambda)\in D^c\times \mathbb{R}$, the upper value and lower value of the game are given respectively by
	$$
	\overline{\mathbb{J}}(x,\lambda):=\inf_{\left(\pi^1, \tau_1\right)\in \Pi_1\times \mathcal{T}}\sup_{\left(\pi^2, \tau_2\right)\in \Pi_2\times \mathcal{T}}\mathbb{J}_\alpha\left(x,\lambda, \pi^1, \tau_1, \pi^2, \tau_2\right) 
	$$
	and
	$$
	\underline{\mathbb{J}}(x,\lambda):=\sup_{\left(\pi^2, \tau_2\right)\in \Pi_2\times \mathcal{T}}\inf_{\left(\pi^1, \tau_1\right)\in \Pi_1\times \mathcal{T}} \mathbb{J}_\alpha\left(x,\lambda, \pi^1, \tau_1, \pi^2, \tau_2\right). 
	$$
	It is clear that  for all $(x,\lambda)\in D^c\times \mathbb{R}$,
	\begin{equation}\label{value inequality}
		\underline{\mathbb{J}}(x,\lambda) 
		\leq \overline{\mathbb{J}}(x,\lambda). 
	\end{equation}
	\noindent{\bf Value of the game:} The game is said to have a value if,  $\overline{\mathbb{J}}(x,\lambda)=\underline{\mathbb{J}}(x,\lambda)$,  for each $(x,\lambda)\in D^c\times \mathbb{R}$. In this case, the common value is denoted by
	$
	\mathbb{J}(x,\lambda).
	$

	\noindent {\bf Saddle point equilibrium:} A strategy $({\pi^*}^1,\tau^*_1)\in \Pi_1\times \mathcal{T}$ is optimal for Player $1$, if for all $(x,\lambda,\pi^2,\tau_2)\in D^c\times \mathbb{R} \times \Pi_2 \times \mathcal{T}$,
	$$
	\underline{\mathbb{J}}(x,\lambda)\geq \mathbb{J}_\alpha\left(\!x,\lambda,{\pi^*}^{1},\tau^*_1,\pi^2,\tau_2\!\right) 
	$$ and
	a strategy $({\pi^*}^2,\tau^*_2)\in \Pi_2\times \mathcal{T}$ is optimal for Player  $2$, if for all $(x,\lambda,\pi^1,\tau_1)\in D^c\times \mathbb{R} \times \Pi_1\times\mathcal{T}$,
	$$
	\overline{\mathbb{J}}(x,\lambda)\leq \mathbb{J}_\alpha\left(\!x,\lambda,\pi^1,\tau_1,{\pi^*}^2,\tau^*_2\!\right). 
	$$
	A pair	$\left(({\pi^*}^1, \tau^{ *}_1),({\pi^*}^2, \tau_2^*)\right)\in (\Pi_1\times \mathcal{T})\times(\Pi_2\times \mathcal{T}) $ is called a saddle point equilibrium, if for all $(x,\lambda)\in D^c\times \mathbb{R}$, $(\pi^1,\tau_1)\in \Pi_1\times \mathcal{T}$ and $(\pi^2,\tau_2)\in \Pi_2\times \mathcal{T}$,
	$$
	\mathbb{J}_\alpha\left(\!x,\lambda, {\pi^*}^1, {\tau}^*_1, \pi^2, \tau_2\!\right)\ \leq \mathbb{J}_\alpha\left(\!x,\lambda, {\pi^*}^1,{\tau}^*_1, {\pi^*}^2, \tau^*_2\!\right)\leq \mathbb{J}_\alpha\left(\!x,\lambda, \pi^1, \tau_1,{\pi^*}^2, \tau^*_2\!\right). 
	$$ This also implies that  $({\pi^*}^1,\tau^*_1)$ is optimal strategy for Player $1$ and $({\pi^*}^2,\tau^*_2)$ is optimal strategy
	for Player $2$.
	
	\subsection{Filtering process}
	To handle the partial observability of the state process, we introduce a sequence of filtering measures describing the conditional distribution of the unobservable state components given the observed history. This construction enables us to reformulate the game as an equivalent completely observable one. Before constructing the updating or filtering measures, we first assume that the semi-Markov kernel $Q$ has a measurable density $q$ with respect to Lebesgue measure on $\mathbb{R}_+$ and some $\sigma$-finite measures $\xi_X$ and $\xi_Y$. Let $\mathcal{P}(E)$ denote the space of all probability measures on the metric space $E$ with the topology of weak convergence. Thus, for $\theta \in \mathbb{R}_+$ and $E \in \mathcal{B}(X\times Y)$ we have 
	$$
	Q(\theta,E\mid x,y,a,b)=\int_E\int_{0}^{\theta}q(t,x^{'},y^{'}\mid x,y,a,b )\xi_X(dx^{'})\xi_Y(dy^{'})dt.
	$$ For any $x,x^{'}\in X, \eta\in \mathcal{P}(Y\times \mathbb{R}), a\in A(x), b\in B(x)$ and $\ell \in \mathbb{R}$ we define a finite measure on $Y\times \mathbb{R}$ by
	$$
	\begin{aligned}
		\beta(G\mid x,\eta,a,b,\theta,x^{'}):&=\int_{{Y\times \mathbb{R}}}\bigg(\int_G q(\theta,x^{'},y^{'}\mid x,y,a,b)\xi_Y(dy^{'})\delta_{\{(\ell-R(x,y,a,b,\theta))e^{\alpha \theta}\}}(d\lambda^{'})\bigg)\eta(dy,d\ell),\\
	\end{aligned}
	$$ for every $G\in \mathcal{B}(Y\times \mathbb{R})$. Normalizing the finite measure $\beta$ we obtain the probability measure on $Y\times \mathbb{R}$ given by
	$$
	\begin{aligned}
		\Psi(G\mid x,\eta,a,b,\theta,x^{'})&:=\frac{\beta(G\mid x,\eta,a,b,\theta,x^{'})}{\int_{{Y\times \mathbb{R}}}\int_{Y\times \mathbb{R}}q(\theta,x^{'},y^{'}\mid x,y,a,b)\xi_Y(dy^{'})\delta_{\{(\ell-R(x,y,a,b,\theta))e^{\alpha \theta}\}}(d\lambda^{'})\eta(dy,d\ell)}\\
		&=\frac{\beta(G\mid x,\eta,a,b,\theta,x^{'})}{\int_{Y}q^{X}(\theta,x^{'}\mid x,y,a,b)\eta^Y(dy)},
	\end{aligned}
	$$ for $G\in \mathcal{B}(Y\times \mathbb{R})$, where $q^{X}(\theta,x^{'}\mid x,y,a,b):=\int_Y q(\theta,x^{'},y^{'}\mid x,y,a,b)\xi_{Y}(dy^{'})$ is the marginal semi-Markov kernel density (with respect to $\xi_Y$) of the next observable state and $\eta^Y(\cdot):=\eta(\cdot,\mathbb{R})$ is the $Y$-marginal distribution of $\eta$. Next, for each $G\in \mathcal{B}(Y\times \mathbb{R})$, $n\geq 0$ and $h_n \in \mathcal{H}_n$, we define a sequence of probability measures on $Y\times \mathbb{R}$ as follows:
	\begin{equation}\label{filtering measure}
		\begin{aligned}
			&\eta_0(G\mid h_0):=(Q_0\otimes \delta_{\lambda})(G)\\
			&\eta_{n+1}(G\mid h_{n+1})=\eta_{n+1}(G\mid h_{n},a_n,b_n,\theta_{n+1},x_{n+1}):=\Psi(G\mid x_{n},\eta_n(\cdot\mid  h_n),a_n,b_n,\theta_{n+1},x_{n+1}).
	\end{aligned} \end{equation}
	The equation \eqref{filtering measure} is called the filtering equation. The measures $\eta_n(\cdot\mid h_n), n\geq 0$ are the associated filtering measures and the $\mathcal{P}(Y\times \mathbb{R})$ valued stochastic process $\{\eta_n(\cdot\mid H_n)\}_{n\geq 0}$ is the corresponding filtering process. Moreover, we define a new stochastic kernel on $\tilde{Q}$ on $\mathbb{R}_+\times X\times \mathcal{P}(Y\times \mathbb{R})$ given $X\times \mathcal{P}(Y\times \mathbb{R})\times A\times B$ as follows:
	\begin{equation}\label{new kernel}
		\begin{aligned}
			&\tilde{Q}(d\theta,dx^{'},d\eta^{'}\mid x,\eta,a,b)\\&:=\delta_{\{\Psi(\cdot\mid  x,\eta,a,b,\theta,x^{'})\}}( d\eta^{'})\int_{{Y\times \mathbb{R}}} Q(d\theta,dx^{'},Y\mid x,y,a,b)\eta(dy,d\ell)\\
			&=\delta_{\{\Psi(\cdot\mid  x,\eta,a,b,\theta,x^{'})\}}( d\eta^{'})Q(d\theta,dx^{'},Y\mid x,\eta^Y,a,b),
		\end{aligned}
	\end{equation}
	where $Q(d\theta,dx^{'},Y\mid x,\mathcal{Y},a,b):=\int_{Y}Q(d\theta,dx^{'},Y\mid x,y,a,b)\mathcal{Y}(dy)$ for any $\mathcal{Y} \in \mathcal{P}(Y)$.
	\begin{prop}\label{proposition 8}
		For each \(n\geq 0\), \((\pi^1,\pi^2)\in\Pi_1\times\Pi_2\), and
		\((x,\lambda)\in X\times\mathbb{R}\), the following assertions hold:
		\begin{enumerate}[topsep=0pt, partopsep=0pt]
			\item[(a)] $\mathbb{P}_{x,\lambda}^{\pi^1,\pi^2}
			\bigl((Y_n,\Lambda_n)\in C \mid H_n\bigr)
			=
			\eta_n(C\mid H_n),
			~\mathbb{P}_{x,\lambda}^{\pi^1,\pi^2}\text{-a.s.}$, for every \(C\in\mathcal{B}(Y\times\mathbb{R})\).
			
			\item[(b)] $	\mathbb{P}_{x,\lambda}^{\pi^1,\pi^2}
			\Bigl(\Theta_{n+1}\in \mathcal{R},
			X_{n+1}\in E,\,
			\eta_{n+1}(\cdot\mid  H_{n+1})\in F
			\,\Big|\, H_n,A_n,B_n
			\Bigr)$\vspace{0.1cm}\\ 
			$=
			\tilde{Q}
			\Bigl(\mathcal{R}\times
			E\times F
			\,\big|\,
			X_n,\eta_n(\cdot\mid  H_n),A_n,B_n
			\Bigr),~
			\mathbb{P}_{x,\lambda}^{\pi^1,\pi^2}\text{-a.s.},
			$
			for every \(\mathcal{R}\in \mathcal{B}(\mathbb{R}_+),~E\in\mathcal{B}(X)\) and
			\(F\in\mathcal{B}(\mathcal{P}(Y\times\mathbb{R}))\). 
	\end{enumerate}	\end{prop}

	\begin{proof}
		$(a)$ This follows from \cite[Theorem 3.1]{wei2025partially}. 
		
		$(b)$ We see that,
		$$
		\begin{aligned}
			& \mathbb{P}_{x, \lambda}^{\pi^1, \pi^2}\left(\Theta_{n+1}\in \mathcal{R}, X_{n+1} \in E, \eta_{n+1}\left(\cdot \mid H_{n+1}\right) \in F \mid H_n, A_n, B_n\right) \\
			= & \mathbb{E}_{x, \lambda}^{\pi^1, \pi^2}\bigg[\mathbb{P}_{x,\lambda}^{\pi^1, \pi^2}\big(\Theta_{n+1}\in \mathcal{R}, X_{n+1} \in E, \eta_{n+1}\left(\cdot \mid H_{n+1}\right) \in F \mid H_n, Y_n, A_n, B_n\big) \bigg| H_n, A_n, B_n\bigg] \\
			= & \int_Y \!\int_Y\! \int_E \!\int_\mathcal{R}\! I_F\!\left(\!\Psi\!\left(\cdot \!\mid\! X_n, \eta_n\left(\cdot \mid H_n\right), A_n, B_n,\theta^{'}, x^{'}\right)\!\right)\! Q(d\theta^{'},dx^{'}, d y^{'}\! \mid \!X_n, y, A_n, B_n) \mathbb{P}_{x, \lambda}^{\pi^1, \pi^2}(Y_n\in d y \mid H_n) \\
			=& \int_{Y\times \mathbb{R}}  \int_E \int_\mathcal{R}\int_F \!\delta_{\left\{\Psi\left(\cdot \mid X_n, \eta_n\left(\cdot \mid H_n\right), A_n, B_n,\theta^{'}, x^{'}\right)\right\}}\!(d\eta^{'}) Q(d\theta^{'},dx^{'}, Y\! \mid \!X_n, y, A_n, B_n) \mathbb{P}_{x, \lambda}^{\pi^1, \pi^2}(d\lambda^{'},d y \mid H_n)\\
			=& \int_Y \int_E \int_\mathcal{R}\int_F \delta_{\left\{\Psi\left(\cdot \mid X_n, \eta_n\left(\cdot \mid H_n\right), A_n, B_n,\theta^{'}, x^{'}\right)\right\}}\!(d\eta^{'}) Q(d\theta^{'},d x^{'}, Y \!\mid \!X_n, y, A_n, B_n) {\eta}_n^Y(dy \mid H_n) \\&\hspace{11.5cm}\mathbb{P}_{x,\lambda}^{\pi^1,\pi^2}\text{-a.s.}\\
			= &~ \tilde{Q}
			\Bigl(\mathcal{R}\times
			E\times F
			\,\big|\,
			X_n,\eta_n(\cdot\mid  H_n),A_n,B_n
			\Bigr) \quad \mathbb{P}_{x,\lambda}^{\pi^1,\pi^2}\text{-a.s}.
		\end{aligned}
		$$  In the above equation the second equality follows from \eqref{property of probability measure 8} and \eqref{filtering measure}, fourth equality follows from Proposition \ref{proposition 8}$(a)$ and the last one is obtained by using \eqref{new kernel}.
	\end{proof}
	\section{Existence of value and saddle point equilibrium} 
	
	We begin this section with a few key results. Using these results, we establish the value of the game and show that the game admits a saddle-point equilibrium. Let $\mathcal{M}_{[0,1]}(E)$ be the set of all measurable functions on the Borel set $E$ whose range lies in $[0,1]$.  
	
	For  $u\in \mathcal{M}_{[0,1]}(D^c\times \mathcal{P}(Y\times\mathbb{R})),x\in D^c,\eta \in  \mathcal{P}(Y\times\mathbb{R}),a\in A(x) $ and $b\in B(x)$ we define 
	\begin{equation}\label{operator 8.2}
		\begin{aligned}
			(S_{a,b}u)(x,\eta):=&\int_{\mathbb{R}_+\times D^c} u(x^{'},\Psi(\cdot\mid x,\eta,a,b,\theta,x^{'})) Q(d\theta,dx^{'},Y\mid x,\eta^Y,a,b)\\+
			& \int_{{Y\times \mathbb{R}}}\int_{\mathbb{R}_+}Q(d\theta,D\times Y\mid x,y,a,b)I_{\{R(x,y,a,b,\theta)\leq 
				\ell\}}\eta(dy,d\ell)\\
			=&\int_{\mathbb{R}_+\times D^c\times \mathcal{P}(Y\times \mathbb{R})} u(x^{'},\eta^{'}) \tilde{Q}(d\theta,dx^{'},d\eta^{'} \mid x,\eta, a,b)\\+& \int_{{Y\times \mathbb{R}}}\int_{\mathbb{R}_+}Q(d\theta,D\times Y\mid x,y,a,b)I_{\{R(x,y,a,b,\theta)\leq 
				\ell\}}\eta(dy,d\ell). 
		\end{aligned}
	\end{equation}
	
	\begin{rem}\label{continuity of S}
		$(i)$ Since $A(x)$ and $B(x)$ are finite, for each $(x,\eta)\in D^c\times \mathcal{P}(Y\times \mathbb{R})$ the map
		$$
		(\mu,\nu)\mapsto \sum_{a\in A(x)}\sum_{b\in B(x)} (S_{a,b}u)(x,\eta)\mu(a)\nu(b),$$ is continuous on $\mathcal{P}(A(x))\times \mathcal{P}(B(x)).$ 
		
		$(ii)$ Observe from \eqref{operator 8.2} that for any $u \in \mathcal{M}_{[0,1]}(D^c\times \mathcal{P}(Y\times\mathbb{R})) $, we have  $0\leq S_{a,b}u\leq 1$ and $S_{a,b}u$ is  measurable on $(D^c\times \mathcal{P}(Y\times\mathbb{R}))$, for each $(a,b)\in A(x)\times B(x)$.
	\end{rem} 
	Now, define another operator on $\mathcal{M}_{[0,1]}(D^c\times \mathcal{P}(Y\times\mathbb{R}))$ given by	
	\begin{equation}\label{operator 8.1}
		\begin{aligned}(Su)(x,\eta):= &\inf_{\mu\in \mathcal{P}(A(x))}\sup_{\nu\in \mathcal{P}(B(x))}\sum_{a\in A(x)}\sum_{b\in B(x)} (S_{a,b}u)(x,\eta)\mu(a)\nu(b)\\
			=&\sup_{\nu\in \mathcal{P}(B(x))}\inf_{\mu\in \mathcal{P}(A(x))}\sum_{a\in A(x)}\sum_{b\in B(x)} (S_{a,b}u)(x,\eta)\mu(a)\nu(b).
		\end{aligned}
	\end{equation} The equality in \eqref{operator 8.1} follows from Remark \ref{continuity of S} and Fan's min-max theorem \cite[Theorem 1(ii)]{fan1953minimax}, as $\mathcal{P}(A(x))$ and $\mathcal{P}(B(x))$ are convex  and the function is linear in $\mu$ and $\nu$.
	
	\begin{rem}\label{measuarbility remark of H}
		Using Remark \ref{continuity of S}$(ii)$, from \eqref{operator 8.1} we can conclude that for any $u \in \mathcal{M}_{[0,1]}(D^c\times \mathcal{P}(Y\times\mathbb{R})) $, we have  $0\leq Su\leq 1$ and $Su$ is  measurable on $(D^c\times \mathcal{P}(Y\times\mathbb{R}))$. Hence $S$ is an operator from $\mathcal{M}_{[0,1]}(D^c\times \mathcal{P}(Y\times\mathbb{R}))$ to itself.
	\end{rem}	
	
	For any $(x,\eta) \in (D^c\times \mathcal{P}(Y\times\mathbb{R}))$ and $k=1,2$, define
	\begin{equation}\label{8 L}
		L_k(x,\eta):=\int_{Y\times\mathbb{R}}I_{[0,\infty)}(\ell-\zeta_k(x,y))\eta(dy,d\ell).
	\end{equation} For $k=1,2$, $L_k(x,\eta)$ denotes the conditional probability under the current filtering measure $\eta$, that if Player $k$ stops immediately, then the terminal reward doesn't exceed the available remaining budget.  
	
	\begin{lem}\label{lemma 8.1}
		Suppose $u \in \mathcal{M}_{[0,1]}(D^c\times \mathcal{P}(Y\times\mathbb{R}))$. Then $u$ satisfies the following inequalities with bilateral constraints
		\begin{equation}\label{8inequalities}
			\begin{aligned}
				\text{(i)}~ & \min\left\{L_1(x,\eta),L_2(x,\eta)\right\} \leq u(x,\eta)\leq L_1(x,\eta), \\[0.6em]
				\text{(ii)}~ & u(x,\eta)-(Su)(x,\eta) \geq 0, \quad \text{if}\ \ L_1(x,\eta)> u(x,\eta), \\[0.6em]
				\text{(iii)}~ & u(x,\eta)-(Su)(x,\eta) \leq 0, \quad \text{if}\ \ L_2(x,\eta)< u(x,\eta),
			\end{aligned}
		\end{equation}
		for each $(x,\eta)\in D^c\times \mathcal{P}(Y\times\mathbb{R})$, if and only if $u$ satisfies 
		\begin{equation}\label{equation}
			u(x,\eta)=\min\left\{L_1(x,\eta),\max\left\{ L_2(x,\eta),(Su)(x,\eta)\right\} \right\}, \quad \forall(x,\eta)\in D^c\times \mathcal{P}(Y\times\mathbb{R}).
		\end{equation}
	\end{lem}
	\begin{proof} Same as the proof of \cite[Lemma 3.1]{zhang2026partialstoppingdiscrete}
	\end{proof}
	
	Next we define another operator $T$ on $\mathcal{M}_{[0,1]}(D^c\times \mathcal{P}(Y\times\mathbb{R}))$ by
	\begin{equation}\label{fixed point operator 8}
		\begin{aligned}
			Tu(x,\eta)=\min\left\{L_1(x,\eta),\max\left\{ L_2(x,\eta),(Su)(x,\eta)\right\} \right\}, 
		\end{aligned}
	\end{equation}
	We now construct a sequence of functions  $\{u_n\}_{n\ge 0}$ recursively as follows:
	\[
	u_0(x,\eta)= L_1(x,\eta), 
	\qquad 
	u_n = T u_{n-1}, \quad n \ge 1.
	\]
	
	\begin{rem}\label{measurability remark of T}
		The mappings 
		$
		~(x,\eta)\mapsto L_k(x,\eta),~ k=1,2,
		$ are measurable and take values in $[0,1]$, and  by Remark \ref{measuarbility remark of H}, $Su \in \mathcal{M}_{[0,1]}(D^c\times \mathcal{P}(Y\times\mathbb{R}))$. Hence, it follows from \eqref{fixed point operator 8} that the operator $T$ is well-defined and maps $\mathcal{M}_{[0,1]}(D^c\times \mathcal{P}(Y\times\mathbb{R}))$ into itself. Consequently, $u_n \in \mathcal{M}_{[0,1]}(D^c\times \mathcal{P}(Y\times\mathbb{R}))$, for every $n\geq 0$.
	\end{rem}

	%
	\begin{lem}\label{lemma 8.2}
		\begin{itemize}
			\item[(a)] 	 The sequence of functions $\left\{u_n\right\}_{n\geq 0} \in \mathcal{M}_{[0,1]}(D^c\times \mathcal{P}(Y\times\mathbb{R}))$ is non-increasing and converges pointwise to a limit function $u^*\in \mathcal{M}_{[0,1]}(D^c\times \mathcal{P}(Y\times\mathbb{R}))$. In particular, for each $(x,\eta)\in D^c\times \mathcal{P}(Y\times\mathbb{R})$,  $$\lim_{n\to \infty}u_n(x,\eta)=u^*(x,\eta).$$ 
			\item[(b)] $u^*$ is a fixed point of the operator $T$, i.e., $T u^*=u^*$.
		\end{itemize}

	\end{lem}
	\begin{proof}
		$(a)$ Since $u_1(x,\eta)=T u_0(x,\eta)=\min\left\{L_1(x,\eta),\max\left\{ L_2(x,\eta),(Su_0)(x,\eta)\right\} \right\},$ it follows that 
		$ u_1(x,\eta)\leq  L_1(x,\eta)=u_0(x,\eta) .$  Now suppose, $u_{n}(x,\eta)\leq u_{n-1}(x,\eta)$. From \eqref{operator 8.1} it is clear that $S$ is monotone in $u$. Thus, we have,
		
		$$
		\begin{aligned}
			u_{n+1}(x,\eta)=T u_n(x,\eta) & =\min\left\{L_1(x,\eta),\max\left\{(Su_n)(x,\eta), L_2(x,\eta)\right\}\right\}\\
			& \leq \min\left\{ L_1(x,\eta),\max\left\{(Su_{n-1})(x,\eta), L_2(x,\eta)\right\}\right\} \\
			& =T u_{n-1}(x,\eta)=u_n(x,\eta) .
		\end{aligned}
		$$

		Thus, by induction the sequence $\left\{u_n\right\}_{n\geq 0}$ is non-increasing. Again from Remark \ref{measurability remark of T} we have $0\leq u_n(x,\eta) \leq 1$, for every $(x,\eta)\!\in \! D^c\times \mathcal{P}(Y\times\mathbb{R}),~n\geq 0$. Hence, for each $(x,\eta)\!\in \! D^c\times \mathcal{P}(Y\times\mathbb{R})$, the limit
		\begin{equation}\label{limit 8}
			u^*(x,\eta):=\lim _{n \rightarrow \infty} u_n(x,\eta),
		\end{equation} exists.
		Moreover, from Remark \ref{measurability remark of T} we see that $u_n \in \mathcal{M}_{[0,1]}(D^c\times \mathcal{P}(Y\times\mathbb{R}))$, for every $n\geq 0$ and $u_n \downarrow u^*$ pointwise. Therefore,  $u^*\in \mathcal{M}_{[0,1]}(D^c\times \mathcal{P}(Y\times\mathbb{R})).$ \\
		(b). As $T$ is monotone in $u$ and $u_n(x,\eta)\geq u^*(x,\eta)$ for all $n\geq 0$, therefore, for each  $(x,\eta)\in D^c\times \mathcal{P}(Y\times\mathbb{R})$ we have,
		$$
		Tu^*(x,\eta) \leq T u_n(x,\eta)=u_{n+1}(x,\eta) .
		$$
		Taking limit $n \rightarrow \infty$ on both sides we get,
		$T u^*(x,\eta) \leq u^*(x,\eta)$ for all $(x,\eta)\in D^c\times \mathcal{P}(Y\times\mathbb{R})$.
		
		Again for fixed $(x,\eta)\in D^c\times \mathcal{P}(Y\times\mathbb{R})$ and for all $n\geq 0$,
		$$	
		\begin{aligned}
			&u_{n+1}(x,\eta)
			= T u_n(x,\eta)\\
			&= \min\left\{ L_1(x,\eta),\max\left\{(Su_n)(x,\eta), L_2(x,\eta)\right\}\right\}\\
			&= \min\Bigl\{  L_1(x,\eta),\;
			\max\Bigl\{
			\inf_{\mu\in \mathcal{P}(A(x))}\sum_{a\in A(x)}\!\sum_{b\in B(x)}\!
			(S_{a,b}u_n)(x,\eta)\mu(a)\nu^*_n(b), L_2(x,\eta)
			\Bigr\}
			\Bigr\}\\
			&\leq  \min\Bigl\{  L_1(x,\eta),\;
			\max\Bigl\{
			\sum_{a\in A(x)}\!\sum_{b\in B(x)}\!
			(S_{a,b}u_n)(x,\eta)\mu(a)\nu^*_n(b), L_2(x,\eta)
			\Bigr\}
			\Bigr\}.
		\end{aligned}
		$$
		The existence of $\nu^*_n \in \mathcal{P}(B(x))$ in the above inequality follows by Remark \ref{continuity of S}$(i)$ and the fact that $\mathcal{P}(B(x))$ is compact. Since $\mathcal{P}(B(x))$ is compact, there exists a sub-sequence
		$\{\nu^*_{n_k}\}_{k\ge1}$ converging to some $\nu^*\in \mathcal{P}(B(x))$.
		The above inequality holds for every $n$, and therefore, in particular, along the sub-sequence
		$\{n_k\}$. Moreover, since $u_n$ converges to $u^*$ pointwise, it follows that for each $(x,\eta)\in(D^c\times \mathcal{P}(Y\times\mathbb{R}))$,
		$u_{n_k+1}(x,\eta)\to u^*(x,\eta)$,  as $k\to\infty$. Applying $n=n_k$ in the above inequality and letting $k\to\infty$, from dominated convergence theorem we get,
		\begin{align*}
			u^*(x,\eta) 
			& \leq \min\Bigl\{  L_1(x,\eta),\;
			\max\Bigl\{
			\sum_{a\in A(x)}\!\sum_{b\in B(x)}\!
			( (S_{a,b}u^*))(x,\eta)\mu(a)\nu^*(b), L_2(x,\eta)
			\Bigr\}
			\Bigr\}\\
			& \leq \min\Bigl\{  L_1(x,\eta),\;
			\max\Bigl\{
			\sup_{\nu\in P(B(x))}\sum_{a\in A(x)}\!\sum_{b\in B(x)}\!
			( (S_{a,b}u^*))(x,\eta)\mu(a)\nu(b), L_2(x,\eta)
			\Bigr\}
			\Bigr\}.
		\end{align*}
		
		Since the above  inequality is true for all $\mu \in \mathcal{P}(A(x))$, we can conclude that for all $(x,\eta)\in D^c\times \mathcal{P}(Y\times\mathbb{R})$,
		$$
		u^*(x,\eta)\leq\min\left\{ L_1(x,\eta),\max\left\{(Su^*)(x,\eta),L_2(x,\eta)\right\}\right\}=Tu^*(x,\eta).$$ Hence,  $u^*$ is a fixed point of the operator $T$.
	\end{proof}
	Before proceeding to the next result, we will introduce some useful notations. 
	\begin{itemize}
		\item For each $(x,\eta)\in D^c\times \mathcal{P}(Y\times\mathbb{R})$, let $\phi^*(x,\eta) \in \mathcal{P}(A(x))$ and $\psi^*( x,\eta) \in \mathcal{P}(B(x))$ be such that
		\begin{equation}\label{8optimal policy for Player $1$}
			\begin{aligned}
				&\inf_{\mu\in \mathcal{P}(A(x))}\sup_{\nu\in \mathcal{P}(B(x))}\sum_{a\in A(x)}\sum_{b\in B(x)}\!
				(S_{a,b}u^*)(x,\eta)\mu(a)\nu(b)\\
				=&\sup_{\nu\in \mathcal{P}(B(x))}\sum_{a\in A(x)}\sum_{b\in B(x)}\!
				(S_{a,b}u^*)(x,\eta)\phi^*(x,\eta)(a)\nu(b).
			\end{aligned}
		\end{equation}
		and
		\begin{equation}\label{8optimal policy for Player $2$}
			\begin{aligned}
				&\sup_{\nu\in \mathcal{P}(B(x))}\inf_{\mu\in \mathcal{P}(A(x))}\sum_{a\in A(x)}\sum_{b\in B(x)}\!
				(S_{a,b}u^*)(x,\eta)\mu(a)\nu(b)\\
				=&\inf_{\mu\in \mathcal{P}(A(x))}\sum_{a\in A(x)}\sum_{b\in B(x)}\!
				(S_{a,b}u^*)(x,\eta)\mu(a)\psi^*(x,\eta)(b).
			\end{aligned}
		\end{equation}
		
		Since $u^*\in \mathcal{M}_{[0,1]}(D^c\times \mathcal{P}(Y\times\mathbb{R}))$, it follows from Remark~\ref{continuity of S} and \cite[Proposition 7.33]{Bertsekas1996stochastic}(see also \cite[Proposition D.5]{hernandezfirst96}), that there exists measurable functions $\phi^*$ and $\psi^*$ satisfying \eqref{8optimal policy for Player $1$} and \eqref{8optimal policy for Player $2$} respectively. 
		
		\item Let ${\pi^{*}}^1:=\left\{{\pi^{*}_n}^1\right\}_{n \geq 0}\in\Pi_1$ and ${\pi^{*}}^2:=\left\{{\pi^{*}_n}^2\right\}_{n \geq 0}\in\Pi_2$ be such that, for all $n\geq 0$
		\begin{equation}\label{8 optimal strategies}
			{\pi^{*}_n}^1(\cdot \mid h_n):={\phi^*}\left(x_n,\eta_n(\cdot\mid h_n)\right)(\cdot)~~ \text { and }~~ {\pi^{*}_n}^2(\cdot \mid h_n):={\psi^*}\left(x_n,\eta_n(\cdot\mid h_n)\right)(\cdot)~.
		\end{equation} 
		
		Define the stopping times associated with Player $1$ and Player $2$, respectively, by,
		\begin{equation}\label{8stopping times}
			\begin{aligned}
				&\tau^*_1:=\inf\left\{n\geq0 :  u^*(X_n,\eta_n(\cdot\mid H_n))\geq L_1(X_n,\eta_n(\cdot\mid H_n))\right\},\\
				&\tau^*_2:=\inf\left\{n\geq0 :  u^*(X_n,\eta_n(\cdot\mid H_n))\leq L_2(X_n,\eta_n(\cdot\mid H_n))\right\},
			\end{aligned}
		\end{equation} where $L_k$, for $k=1,2$, is as in \eqref{8 L}.
	\end{itemize}

	We now present the main theorem of the paper. It establishes that the function $u^*$, defined in Lemma \ref{lemma 8.2} coincides with the value of the game. Moreover,the pair $(({\pi^{*}}^1,\tau^*_1),({\pi^{*}}^2,\tau^*_2))$ is a saddle point equilibrium of the game. To prove the main result, in addition to Assumption \ref{assumption 8.1}, we need another assumption.
	\begin{assumption}\label{assumption 8.2}
		$\mathbb{P}^{\pi^1,\pi^2}_{x,\lambda}(\tau_D<\infty)=1$ for all $(x,\lambda)\in D^c\times \mathbb{R}$ and $(\pi^1,\pi^2)\in \Pi_1 \times \Pi_2$.
	\end{assumption} Thus, Assumption \ref{assumption 8.2} indicates that, no matter what, the game will eventually stop within finite time almost surely. 
	\begin{rem}
		A sufficient condition on the underlying model components ensuring Assumption \ref{assumption 8.2} is that from every state in $D^c$ there is a uniformly positive probability of reaching $D$ in the next transition. More precisely, suppose there exists some constant $\gamma>0$ satisfying $Q(\mathbb{R}_+,D,Y\mid x,y,a,b)\geq \gamma$ for all $(a,b)\in A(x)\times B(x)$ and $x\in D^c$, then Assumption \ref{assumption 8.2} holds. See \cite{bhabak2022semimarkovprobabilitycriteria} for the proof.
	\end{rem}
	
	\begin{thm}\label{Theorem 8}
		Under Assumptions \ref{assumption 8.1} and \ref{assumption 8.2}, the function $u^*$ is the value of the game and the unique fixed point of the operator $T$. Also, the pair $(({\pi^*}^1,{\tau}^*_1), ({\pi^*}^2,\tau^*_2))\in (\Pi_1\times \mathcal{T})\times (\Pi_2\times \mathcal{T}) $ describes a saddle point equilibrium of the game under consideration.
	\end{thm}
	Throughout the remaining section, Assumptions \ref{assumption 8.1} and \ref{assumption 8.2} are in force. The proof of the main theorem relies on the following auxiliary lemmas. 
	\begin{lem}\label{aux lemma}
		Let $\tau_1,\tau_2\in \mathcal{T}$ be any stopping times for Player $1$ and Player $2$. Then for all $x\in D^c$ and $(\pi^1,\pi^2)\in \Pi_1\times \Pi_2$
		\begin{equation}
			\begin{aligned}
				\hspace{-2.5cm}\text{(a)}~&\lim_{m\rightarrow \infty}\sum_{n=0}^{m}\mathbb{E}_{x,\lambda}^{{\pi}^1,\pi^2}\bigg[u^*(X_{n+1},\eta_{n+1}(\cdot\mid H_{n+1}))I_{\{\tau_D>{n+1}\}}I_{\{{\tau}_1\wedge \tau_2\wedge \tau_D=n+1\}}\bigg]\\
				&=\mathbb{E}_{x,\lambda}^{{\pi}^1,\pi^2}\bigg[u^*(X_{{\tau}_1\wedge \tau_2\wedge\tau_D},\eta_{{\tau}_1\wedge \tau_2\wedge \tau_D}(\cdot\mid H_{{\tau}_1\wedge \tau_2\wedge \tau_D}))I_{\{1\leq {\tau}_1\wedge \tau_2<\infty,{\tau}_1\wedge \tau_2\leq \tau_D-1\}}\bigg],\\
				\hspace{-2.5cm}\text{(b)}~&\lim_{m\rightarrow \infty}\sum_{n=0}^{m}\mathbb{E}_{x,\lambda}^{{\pi}^1,\pi^2}\bigg[I_{\{R({X}_n,Y_n,A_n,B_n,\Theta_{n+1})
					\leq \Lambda_n\}}I_{\{\tau_D={n+1}\}}I_{\{{\tau}_1\wedge \tau_2\wedge \tau_D=n+1\}}\bigg]\\
				&=\mathbb{E}_{x,\lambda}^{{\pi}^1,\pi^2}\bigg[I_{\{\sum_{k=0}^{\tau_D-1}e^{-\alpha T_k}R({X}_k,Y_k,A_k,B_k,\Theta_{k+1})
					\leq \lambda\}}I_{\{\tau_D=\tau_1\wedge\tau_2\wedge\tau_D\}}\bigg].
			\end{aligned}
		\end{equation}
	\end{lem}
	\begin{proof}$(a)$~Observe that $I_{\{\tau_D>{n+1}\}}I_{\{{\tau}_1\wedge \tau_2\wedge \tau_D=n+1\}}=I_{\{{\tau}_1\wedge \tau_2=n+1,\tau_D>n+1\}}$ and the events \\ ${\{{\tau}_1\wedge \tau_2=n+1,\tau_D>n+1\}}_{n\geq 0}$ are disjoint. Hence, we have,
		\begin{equation*}\label{8.6'}
			\begin{aligned}
				&\sum_{n=0}^{m}I_{\{\tau_D>{n+1}\}}I_{\{{\tau}_1\wedge \tau_2\wedge \tau_D=n+1\}}=\sum_{n=0}^{m}I_{\{{\tau}_1\wedge \tau_2=n+1,\tau_D>n+1\}}=I_{\{1\leq{\tau}_1\wedge \tau_2\leq m+1,\tau_D>{\tau}_1\wedge \tau_2\}}\\&=I_{\{1\leq{\tau}_1\wedge \tau_2\leq m+1,{\tau}_1\wedge \tau_2\leq \tau_D-1\}}.
			\end{aligned}
		\end{equation*} Then we have,
		\begin{equation}\label{8.8}
			\begin{aligned}
				&\sum_{n=0}^{m}\mathbb{E}_{x,\lambda}^{{\pi}^1,\pi^2}\bigg[u^*({X}_{n+1},\eta_{n+1}(\cdot\mid H_{n+1}))I_{\{\tau_D>{n+1}\}}I_{\{{\tau}_1\wedge \tau_2\wedge \tau_D=n+1\}}\bigg]\\&=\sum_{n=0}^{m}\mathbb{E}_{x,\lambda}^{{\pi}^1,\pi^2}\bigg[u^*({X}_{{\tau}_1\wedge \tau_2\wedge \tau_D},\eta_{{\tau}_1\wedge \tau_2\wedge \tau_D}(\cdot\mid H_{{\tau}_1\wedge \tau_2\wedge \tau_D}))I_{\{\tau_D>{n+1}\}}I_{\{{\tau}_1\wedge \tau_2\wedge \tau_D=n+1\}}\bigg]\\
				&=\mathbb{E}_{x,\lambda}^{{\pi}^1,\pi^2}\bigg[u^*({X}_{{\tau}_1\wedge \tau_2\wedge \tau_D},\eta_{{\tau}_1\wedge \tau_2\wedge \tau_D}(\cdot\mid H_{{\tau}_1\wedge \tau_2\wedge \tau_D}))I_{\{1\leq{\tau}_1\wedge \tau_2\leq m+1,{\tau}_1\wedge \tau_2\leq \tau_D-1\}}\bigg]
			\end{aligned}
		\end{equation}
		
		Taking limit $m \to \infty$ in \eqref{8.8}, by monotone convergence theorem we get,
		\begin{equation}
			\begin{aligned}
				&\lim_{m\rightarrow \infty}\sum_{n=0}^{m}\mathbb{E}_{x,\lambda}^{{\pi}^1,\pi^2}\bigg[u^*(X_{n+1},\eta_{n+1}(\cdot\mid H_{n+1}))I_{\{\tau_D>{n+1}\}}I_{\{{\tau}_1\wedge \tau_2\wedge \tau_D=n+1\}}\bigg]\\
				&=\mathbb{E}_{x,\lambda}^{{\pi}^1,\pi^2}\bigg[u^*(X_{{\tau}_1\wedge \tau_2\wedge \tau_D},\eta_{{\tau}_1\wedge \tau_2\wedge \tau_D}(\cdot\mid H_{{\tau}_1\wedge \tau_2\wedge \tau_D}))I_{\{1\leq {\tau}_1\wedge \tau_2<\infty,{\tau}_1\wedge \tau_2\leq \tau_D-1\}}\bigg].
			\end{aligned}
		\end{equation}
		
		$(b)$~Since the events $\{\tau_D=n+1,{\tau}_1\wedge\tau_2\wedge\tau_D=n+1\}_{n\geq 0}$ are disjoint, we have, 
		\[
		\sum_{n=0}^{m} I_{\{\tau_D=n+1\}}I_{\{{\tau}_1\wedge\tau_2\wedge\tau_D=n+1\}}
		=
		\sum_{n=0}^{m} I_{\{\tau_D=n+1,{\tau}_1\wedge\tau_2\wedge\tau_D=n+1\}}
		=
		I_{\{\bigcup_{n=0}^{m}\{\tau_D=n+1,{\tau}_1\wedge\tau_2\wedge\tau_D=n+1\}\}}.
		\]
		Moreover,
		\[
		\bigcup_{n=0}^{m}\{\tau_D=n+1,{\tau}_1\wedge\tau_2\wedge\tau_D=n+1\}
		=
		\{\tau_D=\tau_1\wedge\tau_2\wedge\tau_D,\ 1\le \tau_1\wedge\tau_2\wedge\tau_D\leq m+1\}.
		\] Hence,
		\begin{equation}\label{8.6''}
			\begin{aligned}
				&\sum_{n=0}^{m}\mathbb{E}_{x,\lambda}^{{\pi}^1,\pi^2}\bigg[I_{\{R({X}_n,Y_n,A_n,B_n,\Theta_{n+1})
					\leq \Lambda_n\}}I_{\{\tau_D={n+1}\}}I_{\{{\tau}_1\wedge \tau_2\wedge \tau_D=n+1\}}\bigg]\\
				&=\sum_{n=0}^{m}\mathbb{E}_{x,\lambda}^{{\pi}^1,\pi^2}\bigg[I_{\{\sum_{k=0}^{({\tau}_1\wedge \tau_2\wedge \tau_D)-1}e^{-\alpha T_k}R({X}_k,Y_k,A_k,B_k,\Theta_{k+1})
					\leq \lambda\}}I_{\{\tau_D={n+1}\}}I_{\{{\tau}_1\wedge \tau_2\wedge \tau_D=n+1\}}\bigg]\\
				&=\mathbb{E}_{x,\lambda}^{{\pi}^1,\pi^2}\bigg[I_{\{\sum_{k=0}^{({\tau}_1\wedge \tau_2\wedge \tau_D)-1}e^{-\alpha T_k}R({X}_k,Y_k,A_k,B_k,\Theta_{k+1})
					\leq \lambda\}}I_{\{\tau_D={\tau}_1\wedge\tau_2\wedge\tau_D,\ 1\le {\tau}_1\wedge\tau_2\wedge\tau_D\le m+1\}}\bigg].
			\end{aligned}
		\end{equation}

		The first equality in the last equation follows by recursively applying the definition of $\Lambda_n$ given in \eqref{stochastic process 8}. Now, taking limit $m\to \infty$ in \eqref{8.6''} and using the fact that $x \in D^c$ (so that $\tau_D\geq 1$ almost surely),  Assumption \ref{assumption 8.2} together with the monotone convergence theorem yields,
		\begin{equation}\label{8.7}
			\begin{aligned}
				&\lim_{m\rightarrow \infty}\sum_{n=0}^{m}\mathbb{E}_{x,\lambda}^{{\pi}^1,\pi^2}\bigg[I_{\{R({X}_n,Y_n,A_n,B_n,\Theta_{n+1})
					\leq \Lambda_n\}}I_{\{\tau_D={n+1}\}}I_{\{{\tau}_1\wedge \tau_2\wedge \tau_D=n+1\}}\bigg]\\
				&=\mathbb{E}_{x,\lambda}^{{\pi}^1,\pi^2}\bigg[I_{\{\sum_{k=0}^{\tau_D-1}e^{-\alpha T_k}R({X}_k,Y_k,A_k,B_k,\Theta_{k+1})
					\leq \lambda\}}I_{\{\tau_D={\tau}_1\wedge\tau_2\wedge\tau_D\}}\bigg].
			\end{aligned}
		\end{equation}
	\end{proof}	
	\begin{lem}\label{8 aux lemma 3.4} Let $({\pi^*}^1,{\pi^*}^2)$ and $(\tau^*_1,\tau^*_2)$ be the stationary strategies and stopping times defined by \eqref{8 optimal strategies} and \eqref{8stopping times}, respectively. Then for each $(x,\lambda)\in X\times \mathbb{R}$ the following assertions hold. 
		\begin{enumerate}
			\item[(a)] For each $\pi^2\in\Pi_2$ and stopping time $\tau_2$ of Player $2$, we have
			\begin{equation}\label{8.14 a}
				\begin{aligned}
					&\mathbb{E}_{x,\lambda}^{{\pi^*}^1,\pi^2}\bigg[I_{\{\sum_{k=0}^{\tau_D-1}e^{-\alpha T_k}R(X_k,Y_k,A_k,B_k,\Theta_{k+1})
						\leq \lambda\}}I_{\{\tau_D={\tau}_1^*\wedge\tau_2\wedge\tau_D\}}\bigg]\\
					&+\mathbb{E}_{x,\lambda}^{{\pi^*}^1,\pi^2}\bigg[u^*(X_{{\tau}_1^*\wedge \tau_2\wedge \tau_D},
					\eta_{{\tau}_1^*\wedge \tau_2\wedge \tau_D}(\cdot\mid H_{{\tau}_1^*\wedge \tau_2\wedge \tau_D}))
					I_{\{{\tau}^*_1\wedge \tau_2\leq \tau_D-1\}}\bigg]\\
					&\leq u^*(x,Q_0 \otimes \delta_\lambda).
				\end{aligned}
			\end{equation}
			
			\item[(b)] For each $\pi^1\in\Pi_1$ and stopping time $\tau_1$ of Player $1$, we have
			\begin{equation}
				\begin{aligned}\label{8.14 b}
					&\mathbb{E}_{x,\lambda}^{\pi^1,{\pi^*}^2}\bigg[I_{\{\sum_{k=0}^{\tau_D-1}e^{-\alpha T_k}R(X_k,Y_k,A_k,B_k,\Theta_{k+1})
						\leq \lambda\}}I_{\{\tau_D=\tau_1\wedge\bar{\tau}_2^*\wedge\tau_D\}}\bigg]\\
					&+\mathbb{E}_{x,\lambda}^{\pi^1,{\pi^*}^2}\bigg[u^*(X_{\tau_1\wedge{\tau}_2^*\wedge\tau_D},
					\eta_{\tau_1\wedge{\tau}_2^*\wedge\tau_D}(\cdot\mid H_{\tau_1\wedge{\tau}_2^*\wedge\tau_D}))
					I_{\{\tau_1\wedge{\tau}_2^*\wedge\tau_D\leq\tau_D-1\}}\bigg]\\
					&\geq u^*(x,Q_0 \otimes \delta_\lambda).
				\end{aligned}
			\end{equation}
		\end{enumerate}
	\end{lem}
	\begin{proof}
		$(a)$~	 We have, 
		\begin{equation}\label{8.3.1}
			\begin{aligned}
				&\mathbb{E}_{x,\lambda}^{{\pi^*}^1,\pi^2}\bigg[u^*(X_{n+1},\eta_{n+1}(\cdot\mid H_{n+1}))I_{\{\tau_D>{n+1}\}}\bigg|\mathcal{F}_n\bigg]\\
				&\hspace{5cm}+\mathbb{E}_{x,\lambda}^{{\pi^*}^1,\pi^2}\bigg[I_{\{R(X_n,Y_n,A_n,B_n,\Theta_{n+1})
					\leq \Lambda_n\}}I_{\{\tau_D={n+1}\}}\bigg|\mathcal{F}_n\bigg]\\
				&=\mathbb{E}_{x,\lambda}^{{\pi^*}^1,\pi^2}\bigg[I_{\{\tau_D>{n}\}}\mathbb{E}_{x,\lambda}^{{\pi^*}^1,\pi^2}\bigg[u^*(X_{n+1},\eta_{n+1}(\cdot\mid H_{n+1}))I_{\{ X_{n+1}\in D^c\}}\big| H_n,A_n,B_n\bigg]\bigg|\mathcal{F}_n\bigg]\\
				&\hspace{2cm}+\mathbb{E}_{x,\lambda}^{{\pi^*}^1,\pi^2}\bigg[I_{\{\tau_D>{n}\}}\mathbb{E}_{x,\lambda}^{{\pi^*}^1,\pi^2}\bigg[I_{\{R(X_n,Y_n,A_n,B_n,\Theta_{n+1})
					\leq \Lambda_n\}}I_{\{ X_{n+1}\in D\}}\big| H_n,A_n,B_n\bigg]\bigg|\mathcal{F}_n\bigg]\\
				&= \Bigg(\sum_{a\in A(X_n)}\sum_{b\in B(X_n)}\int_{{D^c}\times \mathbb{R}_+\times \mathcal{P}(Y\times \mathbb{R})}u^*(x^{'},\eta^{'})\tilde{Q}(d\theta,dx^{'},d\eta^{'}\mid X_n,\eta_{n}(\cdot\mid H_{n}),a,b)\\
				&\hspace{2cm}\phi^*(a\mid X_n,\eta_{n}(\cdot\mid H_{n}))\pi^2_n(b\mid X_0,A_0,B_0,T_1,\ldots,X_{n-1},A_{n-1},B_{n-1},T_n,X_n)\\
				&+\sum_{a\in A(X_n)}\sum_{b\in B(X_n)}\int_{{Y\times \mathbb{R}}}\int_{\mathbb{R}_+}I_{\{R(X_n,y,a,b,\theta)
					\leq \ell\}}Q(d\theta,D\times Y\mid X_n,y,a,b)\eta_n(dy,d\ell\mid H_n)\\
				&\hspace{2cm}\phi^*(a\mid X_n,\eta_{n}(\cdot\mid H_{n}))\pi^2_n(b\mid X_0,A_0,B_0,T_1,\ldots,X_{n-1},A_{n-1},B_{n-1},T_n,X_n)\bigg)I_{\{\tau_D>n\}}\\
				&\leq \sup_{\nu\in \mathcal{P}(B(X_n))}\bigg(\bigg[\sum_{a\in A(X_n)}\sum_{b\in B(X_n)}\int_{{D^c}\times \mathbb{R}_+\times \mathcal{P}(Y\times \mathbb{R})}u^*(x^{'},\eta^{'})\tilde{Q}(d\theta,dx^{'}d\eta^{'}\mid X_n,\eta_{n}(\cdot\mid H_{n}),a,b)\\
				&\hspace{2cm}\phi^*(a\mid X_n,\eta_{n}(\cdot\mid H_{n}))\nu(b)\\
				&+\sum_{a\in A(X_n)}\sum_{b\in B(X_n)}\int_{{Y\times \mathbb{R}}}\int_{\mathbb{R}_+}I_{\{R(X_n,y,a,b,\theta)
					\leq \ell\}}Q(d\theta,D\times Y\mid x,y,a,b)\eta_n(dy,d\ell\mid H_n)\\
				&\hspace{2cm}\phi^*(a\mid X_n,\eta_{n}(\cdot\mid H_{n}))\nu(b)\bigg]\bigg)I_{\{\tau_D>n\}}\\
				& =Su^*(X_n,\eta_n(\cdot\mid H_n))I_{\{\tau_D>n\}}.
			\end{aligned}
		\end{equation} 
		
		In \eqref{8.3.1}, the first term in the second equality is obtained by using Proposition \ref{proposition 8}(b), \eqref{property of probability measure 8} and \eqref{8 optimal strategies}, while the second term follows from Proposition \ref{proposition 8}(a) and \eqref{property of probability measure 8}. Since $\{{\tau}_1^*\wedge \tau_2\wedge \tau_D>n\}$ is $\mathcal{F}_n$ measurable and ${\{{\tau}_1^*\wedge \tau_2\wedge \tau_D>n\}}\subset {\{\tau_D>n\}}$, multipying both side of  \eqref{8.3.1} by $I_{\{{\tau}_1^*\wedge \tau_2\wedge \tau_D>n\}}$ we get, 	
		\begin{equation}\label{8.3.2}
			\begin{aligned}
				&\mathbb{E}_{x,\lambda}^{{\pi^*}^1,\pi^2}\bigg[u^*(X_{n+1},\eta_{n+1}(\cdot\mid H_{n+1}))I_{\{\tau_D>{n+1}\}}I_{\{{\tau}_1^*\wedge \tau_2\wedge \tau_D>n\}}\bigg|\mathcal{F}_n\bigg]\\
				&\hspace{3cm}+\mathbb{E}_{x,\lambda}^{{\pi^*}^1,\pi^2}\bigg[I_{\{R(X_n,Y_n,A_n,B_n,\Theta_{n+1})
					\leq \Lambda_n\}}I_{\{\tau_D={n+1}\}}I_{\{{\tau}_1^*\wedge \tau_2\wedge \tau_D=n+1\}}\bigg| \mathcal{F}_n\bigg]\\
				&\leq Su^*(X_n,\eta_n(\cdot\mid H_n))I_{\{{\tau}_1^*\wedge \tau_2\wedge \tau_D>n\}}. 
			\end{aligned}
		\end{equation}
		From the definition of ${\tau}_1^*$, for each $n<{\tau}_1^*$, we have, $u^*(X_n,\eta_n(\cdot\mid H_n))<L_1(X_n,\eta_n(\cdot\mid H_n))$. Now using Lemma \ref{lemma 8.2}$(b)$ and Lemma \ref{lemma 8.1}(specifically, \eqref{8inequalities}$(ii)$) in \eqref{8.3.2}, and then taking expectation on both sides we obtain,
		\begin{equation}\label{8.1}
			\begin{aligned}
				&\mathbb{E}_{x,\lambda}^{{\pi^*}^1,\pi^2}\bigg[u^*(X_{n+1},\eta_{n+1}(\cdot\mid H_{n+1}))I_{\{\tau_D>{n+1}\}}I_{\{{\tau}_1^*\wedge \tau_2\wedge \tau_D>n\}}\\
				&\hspace{5cm}+I_{\{R(X_n,Y_n,A_n,B_n,\Theta_{n+1})
					\leq \Lambda_n\}}I_{\{\tau_D={n+1}\}}I_{\{{\tau}_1^*\wedge \tau_2\wedge \tau_D=n+1\}}\bigg]\\
				&\leq \mathbb{E}_{x,\lambda}^{{\pi^*}^1,\pi^2}\bigg[u^*(X_{n},\eta_{n}(\cdot\mid H_{n}))I_{\{{\tau}_1^*\wedge \tau_2\wedge \tau_D>n\}}\bigg].
			\end{aligned}
		\end{equation}
		From \eqref{8.1} we get,
		\begin{equation}\label{8.4}
			\begin{aligned}
				&\mathbb{E}_{x,\lambda}^{{\pi^*}^1,\pi^2}\bigg[u^*(X_{n+1},\eta_{n+1}(\cdot\mid H_{n+1}))I_{\{\tau_D>{n+1}\}}I_{\{{\tau}_1^*\wedge \tau_2\wedge \tau_D=n+1\}}\bigg]\\
				&+\mathbb{E}_{x,\lambda}^{{\pi^*}^1,\pi^2}\bigg[u^*(X_{n+1},\eta_{n+1}(\cdot\mid H_{n+1}))I_{\{{\tau}_1^*\wedge \tau_2\wedge \tau_D>{n+1}\}}\bigg]\\&+\mathbb{E}_{x,\lambda}^{{\pi^*}^1,\pi^2}\bigg[I_{\{R(X_n,Y_n,A_n,B_n,\Theta_{n+1})
					\leq \Lambda_n\}}I_{\{\tau_D={n+1}\}}I_{\{{\tau}_1^*\wedge \tau_2\wedge \tau_D=n+1\}}\bigg]\\
				&\leq \mathbb{E}_{x,\lambda}^{{\pi^*}^1,\pi^2}\bigg[u^*(X_{n},\eta_{n}(\cdot\mid H_{n}))I_{\{{\tau}_1^*\wedge \tau_2\wedge \tau_D>n\}}\bigg].
			\end{aligned}
		\end{equation}
		The second term in the left hand side of \eqref{8.4} is obtained from the fact that $\{{\tau}_1^*\wedge \tau_2\wedge \tau_D>{n+1}\}\subset\{\tau_D>{n+1}\} $. Now, put $n=0,1,\ldots,m$ successively in \eqref{8.4}, then summing these $m+1$ inequalities, we get 
		\begin{equation}\label{8.5}
			\begin{aligned}
				&\sum_{n=0}^{m}\mathbb{E}_{x,\lambda}^{{\pi^*}^1,\pi^2}\bigg[u^*(X_{n+1},\eta_{n+1}(\cdot\mid H_{n+1}))I_{\{\tau_D>{n+1}\}}I_{\{{\tau}_1^*\wedge \tau_2\wedge \tau_D=n+1\}}\bigg]\\&+\mathbb{E}_{x,\lambda}^{{\pi^*}^1,\pi^2}\bigg[u^*(X_{m+1},\eta_{m+1}(\cdot\mid H_{m+1}))I_{\{{\tau}_1^*\wedge \tau_2\wedge \tau_D>{m+1}\}}\bigg]\\&+\sum_{n=0}^{m}\mathbb{E}_{x,\lambda}^{{\pi^*}^1,\pi^2}\bigg[I_{\{R(X_n,Y_n,A_n,B_n,\Theta_{n+1})
					\leq \Lambda_n\}}I_{\{\tau_D={n+1}\}}I_{\{{\tau}_1^*\wedge \tau_2\wedge \tau_D=n+1\}}\bigg]\\
				&\leq \mathbb{E}_{x,\lambda}^{{\pi^*}^1,\pi^2}\bigg[u^*(X_{0},\eta_{0}(\cdot\mid H_{0}))I_{\{{\tau}_1^*\wedge \tau_2\wedge \tau_D>0\}}\bigg].
			\end{aligned}
		\end{equation}
		Taking limit $m \to \infty$ in the second term of \eqref{8.5}, using dominated convergence theorem and Assumption \ref{assumption 8.2} we get,
		\begin{equation}\label{8.5'}
			\begin{aligned}
				\lim_{m\rightarrow \infty}\mathbb{E}_{x,\lambda}^{{\pi^*}^1,\pi^2}\bigg[u^*(X_{m+1},\eta_{m+1}(\cdot\mid H_{m+1}))I_{\{{\tau}_1^*\wedge \tau_2\wedge \tau_D>{m+1}\}}\bigg]=0
			\end{aligned}
		\end{equation}
		Then, taking limit $m\to \infty$ in \eqref{8.5}, using Lemma \ref{aux lemma} and \eqref{8.5'} we get,
		\begin{equation}\label{8.11}
			\begin{aligned}
				&\mathbb{E}_{x,\lambda}^{{\pi^*}^1,\pi^2}\bigg[I_{\{\sum_{k=0}^{\tau_D-1}e^{-\alpha T_k}R(X_k,Y_k,A_k,B_k,\Theta_{k+1})
					\leq \lambda\}}I_{\{\tau_D={\tau}_1^*\wedge\tau_2\wedge\tau_D\}}\bigg]\\
				&+\mathbb{E}_{x,\lambda}^{{\pi^*}^1,\pi^2}\bigg[u^*(X_{{\tau}_1^*\wedge \tau_2\wedge \tau_D},\eta_{{\tau}_1^*\wedge \tau_2\wedge \tau_D}(\cdot\mid H_{{\tau}_1^*\wedge \tau_2\wedge \tau_D}))I_{\{1\leq{\tau}_1^*\wedge \tau_2<\infty,{\tau}_1^*\wedge \tau_2\leq \tau_D-1\}}\bigg]\\
				&\leq \mathbb{E}_{x,\lambda}^{{\pi^*}^1,\pi^2}\bigg[u^*(X_{0},\eta_{0}(\cdot\mid H_{0}))I_{\{{\tau}_1^*\wedge \tau_2\wedge \tau_D>0\}}\bigg]
			\end{aligned}.
		\end{equation}
		Also observe that, as $\mathbb{P}_{x,\lambda}^{{\pi^*}^1,\pi^2}(X_0=x)=1$ and $x\in D^c$,  we have $\tau_D\geq 1$. Hence,
		\begin{equation}\label{8.12}
			\begin{aligned}
				&\mathbb{E}_{x,\lambda}^{{\pi^*}^1,\pi^2}\bigg[u^*(X_{0},\eta_{0}(\cdot\mid H_{0}))I_{\{{\tau}_1^*\wedge \tau_2\wedge \tau_D=0\}}\bigg]\\
				&=\mathbb{E}_{x,\lambda}^{{\pi^*}^1,\pi^2}\bigg[u^*(X_{{\tau}_1^*\wedge \tau_2\wedge \tau_D},\eta_{{\tau}_1^*\wedge \tau_2\wedge \tau_D}(\cdot\mid H_{{\tau}_1^*\wedge \tau_2\wedge \tau_D}))I_{\{{\tau}_1^*\wedge \tau_2\leq \tau_D-1,{\tau}_1^*\wedge \tau_2=0\}}\bigg]
			\end{aligned}
		\end{equation}
		
		Adding \eqref{8.11} and \eqref{8.12} we get for all $x\in D^c$,
		\begin{equation}\label{8.13}
			\begin{aligned}
				&\mathbb{E}_{x,\lambda}^{{\pi^*}^1,\pi^2}\bigg[I_{\{\sum_{k=0}^{\tau_D-1}e^{-\alpha T_k}R(X_k,Y_k,A_k,B_k,\Theta_{k+1})
					\leq \lambda\}}I_{\{\tau_D={\tau}_1^*\wedge\tau_2\wedge\tau_D\}}\bigg]\\
				&+\mathbb{E}_{x,\lambda}^{{\pi^*}^1,\pi^2}\bigg[u^*(X_{{\tau}_1^*\wedge \tau_2\wedge \tau_D},\eta_{{\tau}_1^*\wedge \tau_2\wedge \tau_D}(\cdot\mid H_{{\tau}_1^*\wedge \tau_2\wedge \tau_D}))I_{\{0\leq{\tau}_1^*\wedge \tau_2<\infty,{\tau}_1^*\wedge \tau_2\leq \tau_D-1\}}\bigg]\\
				&\leq \mathbb{E}_{x,\lambda}^{{\pi^*}^1,\pi^2}\bigg[u^*(X_{0},\eta_{0}(\cdot\mid H_{0}))\bigg]\\
				&=u^*(x,Q_0 \otimes \delta_\lambda).
			\end{aligned}
		\end{equation}
		Under Assumption \ref{assumption 8.2}, from \eqref{8.13} we finally get,
		\begin{equation*}\label{8.14}
			\begin{aligned}
				&\mathbb{E}_{x,\lambda}^{{\pi^*}^1,\pi^2}\bigg[I_{\{\sum_{k=0}^{\tau_D-1}e^{-\alpha T_k}R(X_k,Y_k,A_k,B_k,\Theta_{k+1})
					\leq \lambda\}}I_{\{\tau_D={\tau}_1^*\wedge\tau_2\wedge\tau_D\}}\bigg]\\
				&+\mathbb{E}_{x,\lambda}^{{\pi^*}^1,\pi^2}\bigg[u^*(X_{{\tau}_1^*\wedge \tau_2\wedge \tau_D},\eta_{{\tau}_1^*\wedge \tau_2\wedge \tau_D}(\cdot\mid H_{{\tau}_1^*\wedge \tau_2\wedge \tau_D}))I_{\{{\tau}_1^*\wedge \tau_2\leq \tau_D-1\}}\bigg]\\
				&\leq u^*(x,Q_0 \otimes \delta_\lambda).
			\end{aligned}
		\end{equation*}
		
		$(b)$~ Similar proof as $(a)$.
	\end{proof}
	
	\begin{lem}\label{lemma 8.3}
		\begin{enumerate}
			\item[(a)]For all $(x,\lambda)\in D^c\times\mathbb{R} $ and $(\pi^2,\tau_2)\in \Pi_2\times \mathcal{T}$,
			\begin{equation}\label{second ultimate 8.1}
				u^*(x,Q_0 \otimes \delta_\lambda)\geq \mathbb{J}_\alpha(x,\lambda,{\pi^*}^1,{\tau}_1^*,\pi^2,\tau_2).
			\end{equation}
			\item[(b)] For all $(x,\lambda)\in D^c\times\mathbb{R} $ and $(\pi^1,\tau_1)\in \Pi_1\times \mathcal{T}$,
			\begin{equation}\label{second ultimate 8.2}
				u^*(x,Q_0 \otimes \delta_\lambda)\leq \mathbb{J}_\alpha(x,\lambda,{\pi}^1,\tau_1,{\pi^*}^2,{\tau}_2^*).
			\end{equation}
		\end{enumerate}
	\end{lem}
	
	\begin{proof}
		\begin{itemize}
			\item[$(a)$]Let, $(x,\lambda)\!\in \!D^c\times\mathbb{R}$ be such that $u^*(x,Q_0 \otimes \delta_\lambda)\!\geq \!L_1(x,Q_0\otimes \delta_\lambda).$ Then,  using $\mathbb{P}_{x,\lambda}^{{\pi^*}^1,\pi^2}(X_0=x)=1$,  it follows from the definition of ${\tau}_1^*$ that ${\tau}_1^*=0$, $\mathbb{P}_{x,\lambda}^{{\pi^*}^1,\pi^2}$-a.s.  Therefore, from \eqref{cost criterion 8} using Proposition \ref{proposition 8}(a) and \eqref{filtering measure} we have,
			\begin{equation}\label{8.3.0}
				\mathbb{J}_\alpha(x,\lambda,{\pi^*}^1,{\tau}_1^*,\pi^2,\tau_2)=\mathbb{P}_{x,\lambda}^{{\pi^*}^1,\pi^2}(\zeta_1(X_0,Y_0)\leq \Lambda_0)=L_1(x,Q_0\otimes \delta_\lambda).
			\end{equation}
			Since we assumed $u^*(x,Q_0 \otimes \delta_\lambda)\!\geq \!L_1(x,Q_0\otimes \delta_\lambda)$ and by Lemma \ref{lemma 8.2}$(b)$ together with Lemma \ref{lemma 8.1}(specifically \eqref{8inequalities}$(i)$)  we have $u^*(x,Q_0\otimes \delta_\lambda)\leq L_1(x,Q_0\otimes \delta_\lambda)$ we have $u^*(x,Q_0\otimes \delta_\lambda)=L_1(x,Q_0\otimes \delta_\lambda)$. Hence, for $(x,\lambda)\!\in \!D^c\times\mathbb{R}$ such that $u^*(x,Q_0 \otimes \delta_\lambda)\!\geq \!L_1(x,Q_0\otimes \delta_\lambda)$, from \eqref{8.3.0} we have,
			\begin{equation}\label{8.3.0.0}
				\mathbb{J}_\alpha(x,\lambda,{\pi^*}^1,{\tau}_1^*,\pi^2,\tau_2)=u^*(x,Q_0 \otimes \delta_\lambda).
			\end{equation}
			
			Now let $(x,\lambda)\!\in \!D^c\times\mathbb{R}$ be such that $u^*(x,Q_0 \otimes \delta_\lambda)<L_1(x,Q_0\otimes \delta_\lambda).$ The first term on the left hand side of \eqref{8.14 a} in Lemma \ref{8 aux lemma 3.4}$(a)$  can be re-written as
			\begin{equation}\label{8.15}
				\begin{aligned}
					&\mathbb{E}_{x,\lambda}^{{\pi^*}^1,\pi^2}\bigg[I_{\{\sum_{k=0}^{\tau_D-1}e^{-\alpha T_k}R(X_k,Y_k,A_k,B_k,\Theta_{k+1})
						\leq \lambda\}}I_{\{\tau_D={\tau}_1^*\wedge\tau_2\wedge\tau_D\}}\bigg]\\
					&=\mathbb{P}_{x,\lambda}^{{\pi^*}^1,\pi^2}\bigg[\sum_{k=0}^{\tau_D-1}e^{-\alpha T_k}R(X_k,Y_k,A_k,B_k,\Theta_{k+1})\leq \lambda,\tau_D\leq {\tau}_1^*, \tau_D\leq \tau_2\bigg].
				\end{aligned}
			\end{equation}
			From the second term on the left hand side of \eqref{8.14 a} in Lemma \ref{8 aux lemma 3.4}$(a)$ we obtain,
			\begin{equation}\label{8.16}
				\begin{aligned}
					&	\mathbb{E}_{x,\lambda}^{{\pi^*}^1,\pi^2}\bigg[u^*(X_{{\tau}_1^*\wedge \tau_2\wedge \tau_D},\eta_{{\tau}_1^*\wedge \tau_2\wedge \tau_D}(\cdot\mid H_{{\tau}_1^*\wedge \tau_2\wedge \tau_D}))I_{\{{\tau}_1^*\wedge \tau_2\leq \tau_D-1\}}\bigg]\\
					&=\mathbb{E}_{x,\lambda}^{{\pi^*}^1,\pi^2}\bigg[u^*(X_{{\tau}_1^*},\eta_{{\tau}_1^*}(\cdot\mid H_{{\tau}_1^*}))I_{\{{\tau}_1^*\leq\tau_2\}}I_{\{{\tau}_1^*\wedge \tau_2\leq \tau_D-1\}}\bigg]\\
					&+\mathbb{E}_{x,\lambda}^{{\pi^*}^1,\pi^2}\bigg[u^*(X_{\tau_2},\eta_{\tau_2}(\cdot\mid H_{ \tau_2}))I_{\{\tau_2<{\tau}_1^*\}}I_{\{{\tau}_1^*\wedge \tau_2\leq \tau_D-1\}}\bigg]
				\end{aligned}
			\end{equation}
			Now observe that 
			\begin{equation}\label{8.17}
				\begin{aligned}
					&\mathbb{E}_{x,\lambda}^{{\pi^*}^1,\pi^2}\bigg[u^*(X_{{\tau}_1^*},\eta_{{\tau}_1^*}(\cdot\mid H_{{\tau}_1^*}))I_{\{{\tau}_1^*\leq\tau_2\}}I_{\{{\tau}_1^*\wedge \tau_2\leq \tau_D-1\}}\bigg]\\
					&=\mathbb{E}_{x,\lambda}^{{\pi^*}^1,\pi^2}\bigg[\sum_{k=0}^{\infty}u^*(X_{k},\eta_{k}(\cdot\mid H_{k}))I_{\{k={\tau}_1^*\leq\tau_2\}}I_{\{{\tau}_1^*\wedge \tau_2\leq \tau_D-1\}}\bigg]\\
					&=\mathbb{E}_{x,\lambda}^{{\pi^*}^1,\pi^2}\bigg[\sum_{k=0}^{\infty}L_1(X_k,\eta_k(\cdot\mid H_k))I_{\{k={\tau}_1^*\leq\tau_2\}}I_{\{{\tau}_1^*\wedge \tau_2\leq \tau_D-1\}}\bigg]\\
					&=\mathbb{E}_{x,\lambda}^{{\pi^*}^1,\pi^2}\bigg[\sum_{k=0}^{\infty}I_{\{k={\tau}_1^*\leq\tau_  2\}}I_{\{{\tau}_1^*\wedge \tau_2\leq \tau_D-1\}}\mathbb{E}_{x,\lambda}^{{\pi^*}^1,\pi^2}\bigg[I_{\{\sum_{i=0}^{k-1}e^{-\alpha T_i}R(X_i,Y_i,A_i,B_i,\Theta_{i+1})+e^{-\alpha T_k}\zeta_1(X_{k},Y_{k})\leq \lambda\}}\bigg| \mathcal{F}_k\bigg]\bigg]\\
					&=\sum_{k=0}^{\infty}\mathbb{E}_{x,\lambda}^{{\pi^*}^1,\pi^2}\bigg[I_{\{\sum_{i=0}^{k-1}e^{-\alpha T_i}R(X_i,Y_i,A_i,B_i,\Theta_{i+1})+e^{-\alpha T_k}\zeta_1(X_{k},Y_{k})\leq \lambda\}}I_{\{k={\tau}_1^*\leq\tau_2\}}I_{\{{\tau}_1^*\wedge \tau_2\leq \tau_D-1\}}\bigg]\\
					&=\sum_{k=0}^{\infty}\mathbb{P}_{x,\lambda}^{{\pi^*}^1,\pi^2}\bigg[\bigg\{\sum_{i=0}^{k-1}e^{-\alpha T_i}R(X_i,Y_i,A_i,B_i,\Theta_{i+1})+e^{-\alpha T_k}\zeta_1(X_{k},Y_{k})\leq \lambda\bigg\}\cap\bigg\{k={\tau}_1^*\leq\tau_2 \bigg\}\\&\hspace{11cm}\cap\bigg\{{\tau}_1^*\wedge \tau_2\leq \tau_D-1 \bigg\}\bigg]\\
					&=\mathbb{P}_{x,\lambda}^{{\pi^*}^1,\pi^2}\bigg[\bigg\{\sum_{i=0}^{{\tau}_1^*-1}e^{-\alpha T_i}R(X_i,Y_i,A_i,B_i,\Theta_{i+1})+e^{-\alpha T_{{\tau}_1^*}}\zeta_1(X_{{\tau}_1^*},Y_{{\tau}_1^*})\leq \lambda\bigg\}\cap\bigg\{{\tau}_1^*\leq\tau_2 \bigg\}\\&\hspace{11cm}\cap\bigg\{{\tau}_1^*\wedge \tau_2\leq \tau_D-1 \bigg\}\bigg]\\
				\end{aligned}
			\end{equation}
			
			The second equality in \eqref{8.17} follows from the fact that in this part we have $\mathbb{P}_{x,\lambda}^{{\pi^*}^1,\pi^2}(\tau^*_1=0)=0$, which implies the terms corresponding to $k=0$ vanishes, along with the definition of ${\tau}_1^*$ together with Lemma \ref{lemma 8.2}$(b)$ and  Lemma \ref{lemma 8.1}(specifically \eqref{8inequalities}$(i)$). The third equality follows from \eqref{8 L}, Proposition \ref{proposition 8}$(a)$ and by recursively applying the definition of $\Lambda_n$ given in \eqref{stochastic process 8}. The final equality follows from the fact that $\{{\tau}_1^*=k \}_{k\geq 0}$ are disjoint events. 
			
			Further note that, 
			\begin{equation}\label{8.18}
				\begin{aligned}
					&\mathbb{E}_{x,\lambda}^{{\pi^*}^1,\pi^2}\bigg[u^*(X_{\tau_2},\eta_{\tau_2}(\cdot\mid H_{ \tau_2}))I_{\{\tau_2<{\tau}_1^*\}}I_{\{{\tau}_1^*\wedge \tau_2\leq \tau_D-1\}}\bigg]
					\\&=\mathbb{E}_{x,\lambda}^{{\pi^*}^1,\pi^2}\bigg[\sum_{k=0}^{\infty}u^*(X_{k},\eta_{k}(\cdot\mid H_{ k}))I_{\{k=\tau_2<{\tau}_1^*\}}I_{\{{\tau}_1^*\wedge \tau_2\leq \tau_D-1\}}\bigg]\\
					&\geq \mathbb{E}_{x,\lambda}^{{\pi^*}^1,\pi^2}\bigg[\sum_{k=0}^{\infty}L_2(X_k,\eta_k(\cdot\mid H_k))I_{\{k=\tau_2<{\tau}_1^*\}}I_{\{{\tau}_1^*\wedge \tau_2\leq \tau_D-1\}}\bigg]\\
					&=\mathbb{P}_{x,\lambda}^{{\pi^*}^1,\pi^2}\bigg[\bigg\{\sum_{i=0}^{\tau_2-1}e^{-\alpha T_i}R(X_i,Y_i,A_i,B_i,\Theta_{i+1})+e^{-\alpha T_{\tau_2}}\zeta_2(X_{\tau_2},Y_{\tau_2})\leq \lambda\bigg\}\cap\bigg\{\tau_2<{\tau}_1^* \bigg\}\\&\hspace{10cm}\cap\bigg\{{\tau}_1^*\wedge \tau_2\leq \tau_D-1 \bigg\}\bigg]\\
				\end{aligned}
			\end{equation}
			The inequality in  \eqref{8.18} follows from using the definition of ${\tau}_1^*$ and Lemma \ref{lemma 8.2}$(b)$. Now using \eqref{8.17} and \eqref{8.18}, from \eqref{8.16} we get,
			\begin{equation}\label{8.19}
				\begin{aligned}
					&	\mathbb{E}_{x,\lambda}^{{\pi^*}^1,\pi^2}\bigg[u^*(X_{{\tau}_1^*\wedge \tau_2\wedge \tau_D},\eta_{{\tau}_1^*\wedge \tau_2\wedge \tau_D}(\cdot\mid H_{{\tau}_1^*\wedge \tau_2\wedge \tau_D}))I_{\{{\tau}_1^*\wedge \tau_2\leq \tau_D-1\}}\bigg]\\
					&\geq \mathbb{P}_{x,\lambda}^{{\pi^*}^1,\pi^2}\bigg[\bigg\{\sum_{i=0}^{{\tau}_1^*-1}e^{-\alpha T_i}R(X_i,Y_i,A_i,B_i,\Theta_{i+1})+e^{-\alpha T_{{\tau}_1^*}}\zeta_1(X_{{\tau}_1^*},Y_{{\tau}_1^*})\leq \lambda\bigg\}\cap\bigg\{{\tau}_1^*\leq\tau_2 \bigg\}\\&\hspace{11cm}\cap\bigg\{{\tau}_1^*\wedge \tau_2\leq \tau_D-1 \bigg\}\bigg]\\
					&+\mathbb{P}_{x,\lambda}^{{\pi^*}^1,\pi^2}\bigg[\bigg\{\sum_{i=0}^{\tau_2-1}e^{-\alpha T_i}R(X_i,Y_i,A_i,B_i,\Theta_{i+1})+e^{-\alpha T_{\tau_2}}\zeta_2(X_{\tau_2},Y_{\tau_2})\leq \lambda\bigg\}\cap\bigg\{\tau_2<{\tau}_1^* \bigg\}\\&\hspace{11cm}\cap\bigg\{{\tau}_1^*\wedge \tau_2\leq \tau_D-1 \bigg\}\bigg]\\
					&=\mathbb{P}_{x,\lambda}^{{\pi^*}^1,\pi^2}\bigg[
					\bigg\{
					\sum_{i=0}^{{\tau}_1^*\wedge \tau_2 -1}
					e^{-\alpha T_i}
					R(X_i,Y_i,A_i,B_i,\Theta_{i+1})
					\\
					&\hspace{1.5cm}+e^{-\alpha T_{{\tau}_1^*\wedge \tau_2}}\bigg(\zeta_1(X_{{\tau}_1^*},Y_{{\tau}_1^*})I_{\{{\tau}_1^*\leq\tau_2\}}
					+\zeta_2(X_{\tau_2},Y_{\tau_2})I_{\{\tau_2<{\tau}_1^*\}}\bigg)
					\le \lambda
					\bigg\}\cap
					\bigg\{
					{\tau}_1^*\wedge \tau_2 \le \tau_D-1
					\bigg\}
					\bigg]
				\end{aligned}
			\end{equation}
			Substituting \eqref{8.19} and \eqref{8.15} in \eqref{8.14 a}, we finally obtain, 
			\begin{equation}\label{8.20}
				\begin{aligned}
					&\mathbb{J}_\alpha(x,\lambda,{\pi^*}^1,{\tau}_1^*,\pi^2,\tau_2)\\
					=&\mathbb{P}_{x,\lambda}^{{\pi^*}^1,\pi^2}\bigg[
					\bigg\{
					\sum_{i=0}^{{\tau}_1^*\wedge \tau_2\wedge \tau_D-1}
					e^{-\alpha T_i}
					R(X_i,Y_i,A_i,B_i,\Theta_{i+1})\\
					&\hspace{1cm}+\bigg(e^{-\alpha T_{{\tau}_1^*\wedge \tau_2}}\left(\zeta_1(X_{{\tau}_1^*},Y_{{\tau}_1^*})I_{\{{\tau}_1^*\leq\tau_2\}}
					+\zeta_2(X_{\tau_2},Y_{\tau_2})I_{\{\tau_2<{\tau}_1^*\}}\right)\bigg)I_{\{{\tau}_1^*\wedge \tau_2 \le \tau_D-1\}}
					\le \lambda
					\bigg\}\bigg]\\
					& \leq u^*(x,Q_0 \otimes \delta_\lambda).
				\end{aligned}
			\end{equation}
			Hence, combining \eqref{8.3.0.0} and \eqref{8.20}, we obtain the desired conclusion.

			\item[$(b)$]	Let, $(x,\lambda)\!\in \!D^c\times\mathbb{R}$ be such that $u^*(x,Q_0 \otimes \delta_\lambda)\leq L_2(x,Q_0 \otimes \delta_\lambda).$ Then,  using $\mathbb{P}_{x,\lambda}^{\pi^1,{\pi^*}^2}(X_0=x)=1$,  it follows from the definition of ${\tau}_2^*$ that ${\tau}_2^*=0$, $\mathbb{P}_{x,\lambda}^{\pi^1,{\pi^*}^2}$-a.s. Then from \eqref{cost criterion 8} we get
			\begin{equation}\label{8.21''}
				\mathbb{J}_\alpha(x,\lambda,\pi^1,\tau_1,{\pi^*}^2,{\tau}_2^*)=\mathbb{P}_{x,\lambda}^{\pi^1,{\pi^*}^2}(\zeta_1(X_0,Y_0)I_{\{\tau_1=0\}}+\zeta_2(X_0,Y_0)I_{\{\tau_1>0\}}\leq \Lambda_0).
			\end{equation} 
			Now $\{\tau_1=0\}\!\in {\mathcal{F}}_0$, where ${\mathcal{F}}_0=\sigma(X_0)$. Since $ X_0=x,~\mathbb{P}_{x,\lambda}^{\pi^1,{\pi^*}^2}$ a.s, hence  $\mathbb{P}_{x,\lambda}^{\pi^1,{\pi^*}^2}\!(\tau_1=0)\in\{0,1\}$. Consequently, by  \eqref{8.21''} either $\mathbb{J}_\alpha(x,\lambda,\pi^1,\tau_1,{\pi^*}^2,{\tau}_2^*)=\mathbb{P}_{x,\lambda}^{\pi^1,{\pi^*}^2}(\zeta_1(X_0,Y_0)\leq \Lambda_0)$ or $\mathbb{J}_\alpha(x,\lambda,\pi^1,\tau_1,{\pi^*}^2,{\tau}_2^*)=\mathbb{P}_{x,\lambda}^{\pi^1,{\pi^*}^2}(\zeta_2(X_0,Y_0)\leq \Lambda_0)$. 
			
			Suppose first that, $\mathbb{J}_\alpha(x,\lambda,\pi^1,\tau_1,{\pi^*}^2,{\tau}_2^*)=\mathbb{P}_{x,\lambda}^{\pi^1,{\pi^*}^2}(\zeta_1(X_0,Y_0)\leq \Lambda_0)$, then using Proposition \ref{proposition 8}$(a)$ and \eqref{filtering measure} along with Lemma \ref{lemma 8.2}$(b)$ and Lemma \ref{lemma 8.1}(in particular, \eqref{8inequalities}$(i)$), we get,
			\begin{equation}\label{8.21'}
				\mathbb{J}_\alpha(x,\lambda,\pi^1,\tau_1,{\pi^*}^2,{\tau}_2^*)=\mathbb{P}_{x,\lambda}^{\pi^1,{\pi^*}^2}(\zeta_1(X_0,Y_0)\leq \Lambda_0)=L_1(x,Q_0 \otimes \delta_\lambda )\geq u^*(x,Q_0 \otimes \delta_\lambda).
			\end{equation} 
			Similarly, if $\mathbb{J}_\alpha(x,\lambda,\pi^1,\tau_1,{\pi^*}^2,{\tau}_2^*)=\mathbb{P}_{x,\lambda}^{\pi^1,{\pi^*}^2}(\zeta_2(X_0,Y_0)\leq \Lambda_0)$, using Proposition \ref{proposition 8}$(a)$, together with the assumption $u^*(x,Q_0 \otimes \delta_\lambda)\leq L_2(x,Q_0 \otimes \delta_\lambda)$, we get,
			\begin{equation}\label{8.21"}
				\mathbb{J}_\alpha(x,\lambda,\pi^1,\tau_1,{\pi^*}^2,{\tau}_2^*)=\mathbb{P}_{x,\lambda}^{\pi^1,{\pi^*}^2}(\zeta_2(X_0,Y_0)\leq \Lambda_0)=L_2(x,Q_0 \otimes \delta_\lambda )\geq u^*(x,Q_0 \otimes \delta_\lambda).
			\end{equation} 
			Therefore, if $(x,\lambda)\!\in \!D^c\times\mathbb{R}$ be such that $u^*(x,Q_0 \otimes \delta_\lambda)\leq L_2(x,Q_0 \otimes \delta_\lambda)$, using \eqref{8.21'} and \eqref{8.21"} from \eqref{8.21''} we have, 
			\begin{equation}\label{8.21}
				\mathbb{J}_\alpha(x,\lambda,\pi^1,\tau_1,{\pi^*}^2,{\tau}_2^*)=\mathbb{P}_{x,\lambda}^{\pi^1,{\pi^*}^2}(\zeta_1(X_0,Y_0)I_{\{\tau_1=0\}}+\zeta_2(X_0,Y_0)I_{\{\tau_1>0\}}\leq \Lambda_0)\geq u^*(x,Q_0 \otimes \delta_\lambda).
			\end{equation} 
			
			Now let $x\in D^c$ and $u^*(x,Q_0 \otimes \delta_\lambda)\!>\!L_2(x,Q_0 \otimes \delta_\lambda).$ Then, from Lemma \ref{8 aux lemma 3.4}$(b)$ we get,
			\begin{equation}\label{8.22}
				\begin{aligned}
					&\mathbb{E}_{x,\lambda}^{\pi^1,{\pi^*}^2}\bigg[I_{\{\sum_{k=0}^{\tau_D-1}e^{-\alpha T_k}R(X_k,Y_k,A_k,B_k,\Theta_{k+1})
						\leq \lambda\}}I_{\{\tau_D=\tau_1\wedge\bar{\tau}_2^*\wedge\tau_D\}}\bigg]\\
					&+\mathbb{E}_{x,\lambda}^{\pi^1,{\pi^*}^2}\bigg[u^*(X_{\tau_1\wedge{\tau}_2^*\wedge  \tau_D},\eta_{\tau_1\wedge{\tau}_2^*\wedge  \tau_D}(\cdot\mid H_{\tau_1\wedge{\tau}_2^*\wedge  \tau_D}))I_{\{\tau_1\wedge{\tau}_2^*\wedge  \tau_D\leq \tau_D-1\}}\bigg]\\
					&\geq u^*(x,Q_0 \otimes \delta_\lambda).
				\end{aligned}
			\end{equation}
			The second term in \eqref{8.22} can be re-written as
			\begin{equation}\label{8.23}
				\begin{aligned}
					&\mathbb{E}_{x,\lambda}^{\pi^1,{\pi^*}^2}\bigg[u^*(X_{\tau_1\wedge{\tau}_2^*\wedge  \tau_D},\eta_{\tau_1\wedge{\tau}_2^*\wedge  \tau_D}(\cdot\mid H_{\tau_1\wedge{\tau}_2^*\wedge  \tau_D}))I_{\{\tau_1\wedge{\tau}_2^*\wedge  \tau_D\leq \tau_D-1\}}\bigg]\\
					&=\mathbb{E}_{x,\lambda}^{\pi^1,{\pi^*}^2}\bigg[u^*(X_{\tau_1},\eta_{\tau_1}(\cdot\mid H_{\tau_1}))I_{\{ \tau_1\leq{\tau}_2^* \}}I_{\{\tau_1\wedge{\tau}_2^*\leq \tau_D-1\}}\bigg]\\
					&+\mathbb{E}_{x,\lambda}^{\pi^1,{\pi^*}^2}\bigg[u^*(X_{{\tau}_2^*},\eta_{{\tau}_2^*}(\cdot\mid H_{{\tau}_2^*}))I_{\{ {\tau}_2^*<\tau_1 \}}I_{\{\tau_1\wedge{\tau}_2^*\leq \tau_D-1\}}\bigg].
				\end{aligned}
			\end{equation}
			
			By Lemmas \ref{lemma 8.2} and \ref{lemma 8.1} we have $u^*(X_{k},\eta_{k}(\cdot\mid H_{k}))\leq L_1(X_k,\eta_k(\cdot\mid H_k) $, for any $k\geq 0$. Also from the definition of ${\tau}_2^*$, for any $k={\tau}_2^*$, we have $u^*(X_{k},\eta_{k}(\cdot\mid H_{k}))\leq L_2(X_k,\eta_k(\cdot\mid H_k).$ Using these two facts, mimicking the steps in \eqref{8.17} and \eqref{8.18}, from \eqref{8.23} alongside \eqref{8.21} the lemma follows.
		\end{itemize}
	\end{proof}
	
	\begin{proof}[\textbf{Proof of Theorem \ref{Theorem 8}}]
		From Lemma \ref{lemma 8.3}$(a)$, for each $(x,\lambda)\in D^c\times \mathbb{R}$, we have, 
		\begin{equation}\label{ultimate 8.1}
			\begin{aligned}	u^*(x,Q_0 \otimes \delta_\lambda)\geq \sup_{\left(\pi^2, \tau_2\right)\in \Pi_2\times \mathcal{T}}\mathbb{J}_\alpha&(x,\lambda,{\pi^*}^1,{\tau}_1^*,\pi^2,\tau_2)\\
				&\geq \inf_{\left(\pi^1, \tau_1\right)\in \Pi_1\times \mathcal{T}}\sup_{\left(\pi^2, \tau_2\right)\in \Pi_2\times \mathcal{T}}\mathbb{J}_\alpha(x,\lambda,{\pi}^1,{\tau}_1,\pi^2,\tau_2)=\overline{\mathbb{J}}(\lambda,x).
		\end{aligned}\end{equation} 
		
		Also from Lemma \ref{lemma 8.3}$(b)$, for each $(x,\lambda)\in D^c\times \mathbb{R}$, we get,
		\begin{equation}\label{ultimate 8.2}
			\begin{aligned}u^*(x,Q_0 \otimes \delta_\lambda)\leq \inf_{\left(\pi^1, \tau_1\right)\in \Pi_1\times \mathcal{T}}\mathbb{J}_\alpha&(x,\lambda,{\pi}^1,{\tau}_1,{\pi^*}^2,{\tau}_2^*)\\
				&\leq \sup_{\left(\pi^2, \tau_2\right)\in \Pi_2\times \mathcal{T}}\inf_{\left(\pi^1, \tau_1\right)\in \Pi_1\times \mathcal{T}}\mathbb{J}_\alpha(x,\lambda,{\pi}^1,{\tau}_1,\pi^2,\tau_2)=\underline{\mathbb{J}}(\lambda,x).
			\end{aligned}
		\end{equation}
		
		Using \eqref{value inequality}, from \eqref{ultimate 8.1} and \eqref{ultimate 8.2} it follows that, for each $(x,\lambda)\in D^c\times \mathbb{R}$ ,  $u^*(x,Q_0 \otimes \delta_\lambda)=\overline{\mathbb{J}}(\lambda,x)=\underline{\mathbb{J}}(\lambda,x).$ Hence,  the value of the game exists and is given by $u^*$. This also implies that $u^*$ is the unique fixed point of $T$. Since \eqref{second ultimate 8.1} and \eqref{second ultimate 8.2} hold for any $(\pi^1,\tau_1)\in (\Pi_1\times \mathcal{T})$ and $(\pi^2,\tau_2)\in (\Pi_2\times \mathcal{T})$ respectively, therefore putting $({\pi^*}^1,\tau^*_1)$ in \eqref{second ultimate 8.1} and $({\pi^*}^2,\tau^*_2)$ in \eqref{second ultimate 8.2} we get
		\begin{equation}\label{8 value}
			u^*(x,Q_0 \otimes \delta_\lambda)= \mathbb{J}_{\alpha}(x,\lambda,{\pi^*}^1,{\tau}^*_1, {\pi^*}^2,{\tau}_2^*).
		\end{equation}
		
		Using \eqref{8 value}, from \eqref{second ultimate 8.1} and \eqref{second ultimate 8.2} we can conclude that $(({\pi^*}^1,{\tau}^*_1), ({\pi^*}^2,\tau^*_2))\in (\Pi_1\times \mathcal{T})\times (\Pi_2\times \mathcal{T}) $ forms a saddle point equilibrium.

	\end{proof}


	\bibliographystyle{plain}
	\bibliography{semi_Markov_partially_observable_stopping}

\begin{thebibliography}{10}

\bibitem{baauerle2017partially}
Nicole B{\"a}auerle and Ulrich Rieder.
\newblock Partially observable risk-sensitive {M}arkov decision processes.
\newblock {\em Mathematics of Operations Research}, 42(4):1180--1196, 2017.

\bibitem{Bertsekas1996stochastic}
Dimitri Bertsekas and Steven~E Shreve.
\newblock {\em Stochastic optimal control: the discrete-time case}, volume~5.
\newblock Athena Scientific, 1996.

\bibitem{bhabak2022semimarkovprobabilitycriteria}
Arnab Bhabak, Chandan Pal, and Subhamay Saha.
\newblock Zero-sum semi-{M}arkov games with a probability criterion.
\newblock {\em Stochastics}, 94(3):415--431, 2022.

\bibitem{bhabak2024partially}
Arnab Bhabak and Subhamay Saha.
\newblock Partially observable discrete-time discounted {M}arkov games with general utility.
\newblock {\em Operations Research Letters}, 54:107113, 2024.

\bibitem{boda2004stochastic}
Kang Boda, Jerzy~A Filar, Yuanlie Lin, and Lieneke Spanjers.
\newblock Stochastic target hitting time and the problem of early retirement.
\newblock {\em IEEE Transactions on Automatic Control}, 49(3):409--419, 2004.

\bibitem{bose2026semi}
Bivakar Bose, Chandan Pal, and Subhamay Saha.
\newblock Semi-{M}arkov stopping games with risk-sensitive discounted cost criterion.
\newblock {\em Dynamic Games and Applications}, pages 1--20, 2026.

\bibitem{chen2025zero}
Fang Chen, Xianping Guo, Xin Guo, and Qing-Shan Jia.
\newblock Zero-sum discounted semi-{M}arkov games with stopping and control.
\newblock {\em Discrete and Continuous Dynamical Systems-S}, 18(10):2810--2827, 2025.

\bibitem{deep2023partially}
Akash Deep, Shiyu Zhou, Dharmaraj Veeramani, and Yong Chen.
\newblock Partially observable {M}arkov decision process-based optimal maintenance planning with time-dependent observations.
\newblock {\em European journal of operational research}, 311(2):533--544, 2023.

\bibitem{fan1953minimax}
Ky~Fan.
\newblock Minimax theorems.
\newblock {\em Proceedings of the National Academy of Sciences}, 39(1):42--47, 1953.

\bibitem{feinberg2016partially}
Eugene~A Feinberg, Pavlo~O Kasyanov, and Michael~Z Zgurovsky.
\newblock Partially observable total-cost {M}arkov decision processes with weakly continuous transition probabilities.
\newblock {\em Mathematics of Operations Research}, 41(2):656--681, 2016.

\bibitem{ghosh2004partialzero}
Mrinal~K Ghosh, D~McDonald, and Sagnik Sinha.
\newblock Zero-sum stochastic games with partial information.
\newblock {\em Journal of optimization theory and applications}, 121(1):99--118, 2004.

\bibitem{ghosh2007zero}
Mrinal~K Ghosh and KS~Mallikarjuna Rao.
\newblock Zero-sum stochastic games with stopping and control.
\newblock {\em Operations research letters}, 35(6):799--804, 2007.

\bibitem{hernandezfirst96}
On{\'e}simo Hern{\'a}ndez-Lerma and Jean~B Lasserre.
\newblock {\em Discrete-time {M}arkov control processes: basic optimality criteria}, volume~30.
\newblock Springer Science \& Business Media, 1996.

\bibitem{huang2015firstpassageriskprobabilitysemimarkov}
XiangXiang Huang, XiaoLong Zou, and XianPing Guo.
\newblock A minimization problem of the risk probability in first passage semi-{M}arkov decision processes with loss rates.
\newblock {\em Science China Mathematics}, 58(9):1923--1938, 2015.

\bibitem{huang2011}
Yong-hui Huang and Guo Xian-ping.
\newblock First passage models for denumerable semi-{M}arkov decision processes with nonnegative discounted costs.
\newblock {\em Acta Mathematicae Applicatae Sinica, English Series}, 27(2):177--190, 2011.

\bibitem{huang2013minimum}
Yonghui Huang, Xianping Guo, and Zhongfei Li.
\newblock Minimum risk probability for finite horizon semi-{M}arkov decision processes.
\newblock {\em Journal of Mathematical Analysis and Applications}, 402(1):378--391, 2013.

\bibitem{huang2011performance}
Yonghui Huang, Xianping Guo, and Xinyuan Song.
\newblock Performance analysis for controlled semi-{M}arkov systems with application to maintenance.
\newblock {\em Journal of optimization theory and applications}, 150(2):395--415, 2011.

\bibitem{mamabolo2004maintenance}
Reshoketswe~M Mamabolo and Frank~E Beichelt.
\newblock Maintenance policies with minimal repair.
\newblock {\em Stochastics and Quality Control}, 19(2):143--166, 2004.

\bibitem{ohtsubo2002optimal}
Yoshio Ohtsubo and Kenji Toyonaga.
\newblock Optimal policy for minimizing risk models in {M}arkov decision processes.
\newblock {\em Journal of mathematical analysis and applications}, 271(1):66--81, 2002.

\bibitem{pal2020continuous}
Chandan Pal and Subhamay Saha.
\newblock Continuous-time zero-sum stochastic game with stopping and control.
\newblock {\em Operations Research Letters}, 48(6):715--719, 2020.

\bibitem{parthasarathy2020stochastic}
Thiruvenkatachari Parthasarathy, Sujatha Babu, et~al.
\newblock {\em Stochastic games and related concepts}.
\newblock Springer, 2020.

\bibitem{puterman2014markov}
Martin~L Puterman.
\newblock {\em Markov decision processes: discrete stochastic dynamic programming}.
\newblock John Wiley \& Sons, 2014.

\bibitem{rendleman1987assessing}
Richard~J Rendleman~Jr and Richard~W McEnally.
\newblock Assessing the costs of portfolio insurance.
\newblock {\em Financial Analysts Journal}, 43(3):27--37, 1987.

\bibitem{wei2025partially}
Qingda Wei and Xian Chen.
\newblock Partially observable discrete-time stochastic games under risk probability criterion.
\newblock {\em International Journal of Game Theory}, 54(2):28, 2025.

\bibitem{wei2026firstpassageriskprobability}
Qingda Wei and Xian Chen.
\newblock First passage risk probability optimality for partially observable discrete-time zero-sum games.
\newblock {\em Stochastics}, pages 1--24, 2026.

\bibitem{wen2025riskkybernetika}
Xian Wen, Jinhua Cui, and Haifeng Huo.
\newblock The risk probability optimal problem for infinite discounted semi-{M}arkov decision processes.
\newblock {\em Kybernetika}, 61(4):447--466, 2025.

\bibitem{wen2022optimal}
Xin Wen, Xianping Guo, and Li~Xia.
\newblock Optimal dividend problems with a risk probability criterion.
\newblock {\em Naval Research Logistics (NRL)}, 69(3):421--430, 2022.

\bibitem{wen2024semimarkovpartially}
Xin Wen, Xianping Guo, and Li~Xia.
\newblock Finite horizon partially observable semi-{M}arkov decision processes under risk probability criteria.
\newblock {\em Operations Research Letters}, 57:107187, 2024.

\bibitem{zhang2024discrete}
Wenzhao Zhang and Congying Liu.
\newblock Discrete-time stopping games with risk-sensitive discounted cost criterion.
\newblock {\em Mathematical Methods of Operations Research}, 100(2):437--466, 2024.

\bibitem{zhang2026partialstoppingdiscrete}
Wenzhao Zhang, Jinwen Xu, and Qing Li.
\newblock Discrete-time partially observable stopping games with a risk probability criterion.
\newblock {\em Statistics \& Probability Letters}, page 110727, 2026.

\end{thebibliography}
	
\end{document}